\documentclass[12pt]{amsart}
\usepackage{amsmath, amsfonts, amssymb, amsthm}
\usepackage{mathtools}
\usepackage{graphicx}
\usepackage{caption}
\usepackage{subcaption}
\usepackage{multirow}
\usepackage{url}

\newtheorem{thm}{Theorem}[section]
\newtheorem{conj}[thm]{Conjecture}

\makeatother

\usepackage[utf8]{inputenc}

\title{Data-Driven Perspectives on Knot Invariants }
\author{Pawe{\l} D{\l}otko}
\thanks{PD and DG acknowledges support by Dioscuri program initiated by the Max Planck Society, jointly managed with the National Science Centre (Poland), and mutually funded by the Polish Ministry of Science and Higher Education and the German Federal Ministry of Education and Research}
\address{Dioscuri Centre in Topological Data Analysis\\
	Mathematical Institute, Polish Academy of Sciences\\
	Warsaw, Poland}
\email{pdlotko@impan.pl}

\author{Davide Gurnari}
\address{Dioscuri Centre in Topological Data Analysis\\
	Mathematical Institute, Polish Academy of Sciences\\
	Warsaw, Poland}
\email{dgurnari@impan.pl}

\author{Radmila Sazdanovic}
\thanks{RS was partially supported by NSF grant DMS-185470.}
\address{Department of Mathematics\\
	North Carolina State University\\
	Raleigh, NC}
\email{rsazdanovic@math.ncsu.edu}
\date{\today}

\begin{document}
	
	\begin{abstract}
		Data science offers a powerful tool to understand objects in multiple sciences. In this paper we utilize concept of data science, most notably topological data analysis, to extend our understanding of knot theory. This approach provides a way to extend  mathematical exposition of various invariants of knots towards understanding their relations in  statistical and cumulative way. Paper included examples illustrating how topological data analysis can illuminate  structure and relations between knot invariants, state new hypothesis, and gain new insides into long standing conjectures. 
	\end{abstract}
	
	\maketitle
	%\tableofcontents

	%%%%%%%%%%%%%%%%%%%%%%%%%%%%%
	%\input{intro}
	\section{Introduction}\label{sec:intro}
	
	Throughout the 20th century, and notably accelerating over the past four decades, knot theory has undergone rapid development, catalyzed by milestones such as the invention of the Jones polynomial \cite{jones1985polynomial,jones1987Hecke, jones1997polynomial}, its generalizations such as the HOMFLYPT \cite{freyd1985new, przytycki1987invariants} polynomial, and later its categorification, Khovanov link homology  \cite{khovanov1999categorification, bar2002khovanov}. 
	Theoretical advancements, combined with progress in algorithms and technology, have enabled computation of these invariants for numerous knots, leading to the classification of prime  knots with up to 20 crossings \cite{conway1970enumeration, hoste1998first, hoste2005enumeration,regina, burton2020next, thistlethwaite57enumeration} yielding a total of 1, 847, 318, 507 prime hyperbolic knots.
	
	Despite the vast availability of data and the ability to generate more on demand, data science methods applied to knots and their invariants have encountered significant challenges in attempt to unknot hard knots \cite{applebaum2024unknotting, Gukov_2021},  predict knot invariants \cite{hughes2016neural, craven2023learning}, find relations between them and addressing longstanding open problems \cite{davies2021advancing,jejjala2019deep, levitt2019big, craven2024illuminating, lacabanne2024big},  or use for knot identification \cite{sleiman2024geometric}. However, the theorem of Ernst and Summers \cite{ernst1987growth} proves that the number of knots grows exponentially with the number of crossings making it impossible to obtain general answers using purely computational tools.
	
	In this project, we largely move away from black-box data science methods and instead utilize explainable tools from standard statistics and topological data visualization, such as Mapper~\cite{singh2007topological} and Ball Mapper~\cite{ball_mapper} \cite{dlotko2024mapper}.These approaches, rather than providing answers that are valid only for the specific collections of knots considered, allow us to formulate hypotheses that can be explored theoretically and potentially extended for infinite collection of knots. In this way, we bridge data science-seeking inspiration-with mathematical tools capable of handling an infinite amount of data, which grows exponentially with the number of crossings.
	
	The paper is organized as follows: In Section~\ref{sec:data}, we provide a brief introduction to the datasets considered and the methods used in the subsequent analysis. Section~\ref{sec:intuition} explores the trends in the discriminative power of various knot invariants as the number of crossings increases. Section~\ref{sec:BM} presents topological visualization tools applied to knot data. Section~\ref{sec:BMvsPCA} compares topological visualization with standard principal component analysis. Section~\ref{sec:s_vs_sgnature} discusses data science insights into the relationship between the knot signature and Rasmussen s-invariant and provides interpretations of Ball Mapper graphs. Finally, Section \ref{random} provides a small example of data consisting of the Jones polynomials of random knots, to complement the data analysis on knots filtered by the crossing number.

	%%%%%%%%%%%%%%%%%%%%%%%%%%%%%
	% \section{data }\label{sec:data}
	\section {Data and methods} \label{sec:data}
	
	\subsection{Knot invariants as data}
	
	In this work we step away from considering knot as a geometric object and analyze its numerical characteristics and polynomials being the invariants of the knot. The invariants considered here were computed by J.S. Levitt, D. Shuetz, and D. Gurnari and include Khovanov homology, the Jones, Alexander, and HOMFLYPT polynomials, the minimal crossing number, the alternating property, the signature, and the Rasmussen invariant. For further details, references, as well as access to the datasets, please refer to~\cite{datasets}.
	
	To utilize data science tools, the knot polynomials are converted into point clouds. To achieve this, knot polynomials are represented as vectors of their coefficients, as detailed in \cite{levitt2019big, dlotko2024mapper}. When large collection of polynomials are considered, the obtained vectors are padded with zeros to ensure that all of them have the same length, which is determined by the overall minimum and maximum exponents of the specific polynomial over the considered collection of knots. Additionally, we ensure that the coefficients corresponding to a specific power will occupy the same position in the vector for all the considered knots.
	In case of a two-variable polynomial, such as the HOMFLYPT or Khovanov homology polynomial, the same procedure is used to obtain a coefficient matrix, which is subsequently flattened.
	The resulting coefficient vectors have $17$ entries for Alexander, $51$ for Jones, $152$ for HOMFLYPT and $726$ for the Khovanov homology polynomial. Knots are labeled using the knot ids from the Regina software~\cite{regina}.
	
	\subsection{Random knot data}\label{randomdata}
	In addition to organizing these datasets by crossing number, we augment with a collection of approximately 10,000 random knots, referred to in the rest of the paper as random data set. They are  generated from random polygons of chain length 500, as introduced in \cite{cantarella2024faster}. These knots were simplified using \emph{REAPR: Fast Knot Diagram Simplification and Analysis} by Jason Cantarella, Henrik Schumacher, and Clayton Shonkwiler, with Jones polynomials computed via Regina \cite{regina}. Approximately 1\% of the knots had to be discarded due to aborted Jones polynomial computations.
	%}

\subsection{Mapper-like algorithms }
Mapper and ball mapper are topological data analysis techniques used
to visualize high dimensional complicated data. Since detailed description is provided  in~\cite{dlotko2024mapper} we provide just a brief introduction to the tools focusing on ball mapper and its extensions.
Given a collection of points $X \subset \mathbb{R}^n$ and an $\epsilon
> 0$, an \emph{epsilon net} of $X$ is a subset $Y \subset X$ having
the property that for every $x \in X$ there exist $y \in Y$ in a
distance no greater than $\epsilon$ form $x$. As a consequence, $X
\subset \bigcup_{y \in Y} B(y,\epsilon)$, where $B(y,\epsilon)$ denotes a
ball of a radius $\epsilon$ centered in $y$. Ball mapper is then
defined as an abstract graph whose vertices correspond to points in
$Y$. With a little abuse of notation, we will set $V = Y$, even though
the ball mapper is not an embedded graph. The edges of ball mapper
graph are placed between vertices $y_1, y_2 \in Y$ if there exist a
point $x \in X$ with the property that $d(x,y_1) \leq \epsilon$ and
$d(x,y_2) \leq \epsilon$. In this case, point $x$ is covered mutually
by balls centered in $y_1$ and $y_2$. Such a graph can be visualized
using common libraries for graph visualization. 

When the points of $X$ are additionally equipped with a function $f: X \rightarrow \mathbb{R}$, the values of the function can be lifted to vertices of the ball mapper. More precisely, the value of the vertex $y \in Y$ is set to an aggregation, typically an average, of values of $f$ for all points in $X \cap B(y,\epsilon)$. This aggregated value of $f$  can be presented with an appropriate color scale as a color of a vertex $y$. This way, we obtain a \emph{coloring function} induced by $f$. This way, ball mapper can serve as a way to visualize the point cloud $X \subset \mathbb{R}^n$ as well as a function $f: X \rightarrow \mathbb{R}^n$. 

The visualization of scalar-valued functions can be, as discussed in~\cite{datasets}, extended to visualize vector valued functions. Let us be more precise. Suppose we are given $K$, a collection of knots. Let $f(K)$ and $g(K)$ denote two types of knot invariants in form of potentially high dimensional point clouds. They may, for instance, be the Alexander and Jones polynomials of knots in $K$ expressed in an appropriate basis. Given that both $f(K)$ and $g(K)$ are originated from $K$, they can be
attached by a relation. Namely, $x \in f(K)$ is in relation with $y
\in g(K)$ if there is a knot $k \in K$ such that $f(k) = x$ and $g(k)
= y$.

This point-based relation can be extended to ball mappers $G_{f(K)}$ and $G_{g(K)}$, which are constructed on $f(K)$ and $g(K)$ for appropriately chosen parameters $\epsilon_f$ and $\epsilon_g$ (which may differ). Given a point $y \in f(K)$ in the $\epsilon_f$-grid, we consider all points within $B(y, \epsilon_f)$. Those are all points covered by the vertex $y$ of the $G_{f(K)}$ graph. Given these points, we define a function that takes the value $1$ for all points in $B(y, \epsilon_f)$ and $0$ for all other points. The coloring function induced by it can be applied to $G_{g(K)}$ to point out the location of the vertex $y$ of the ball mapper $G_{f(K)}$ in the ball mapper $G_{g(K)}$. In other words, this approach highlights where the knots grouped in $B(y, \epsilon_f)$ appear within $G_{g(K)}$. 
Note that this mapping can be extended to an arbitrary set of vertices in $G_{f(K)}$. Also, the roles of $f$ and $g$ can be interchanged. Such collection of visualizations provides valuable intuition about the relationship between the two mappings, $f$ and $g$, when defined on the same dataset $K$ .

%\input{dataMethods}

%%%%%%%%%%%%%%%%%%%%%%%%%%%%%
% \section{Intuition- data base questions}\label{sec:intuition}
%\input{intuition}

\section{Statistical insights into the discriminative power of polynomial knot invariants}
\label{sec:intuition}

Every newly introduced knot invariant raises the same fundamental questions: What information does it encode about knots? Is it a complete invariant? Which knots can it distinguish with it? How does it relate to other known invariants? What is its computational complexity? Even the questions of unknot recognition, if there exist a non-trivial knot which has the same value of that invariant as the unknot, is still open for the Jones polynomial and was settled for Khovanov homology \cite{kronheimer2011khovanov} only in 2010.

\begin{figure}[ht!]
	\centering
	\includegraphics[width=0.8\textwidth]{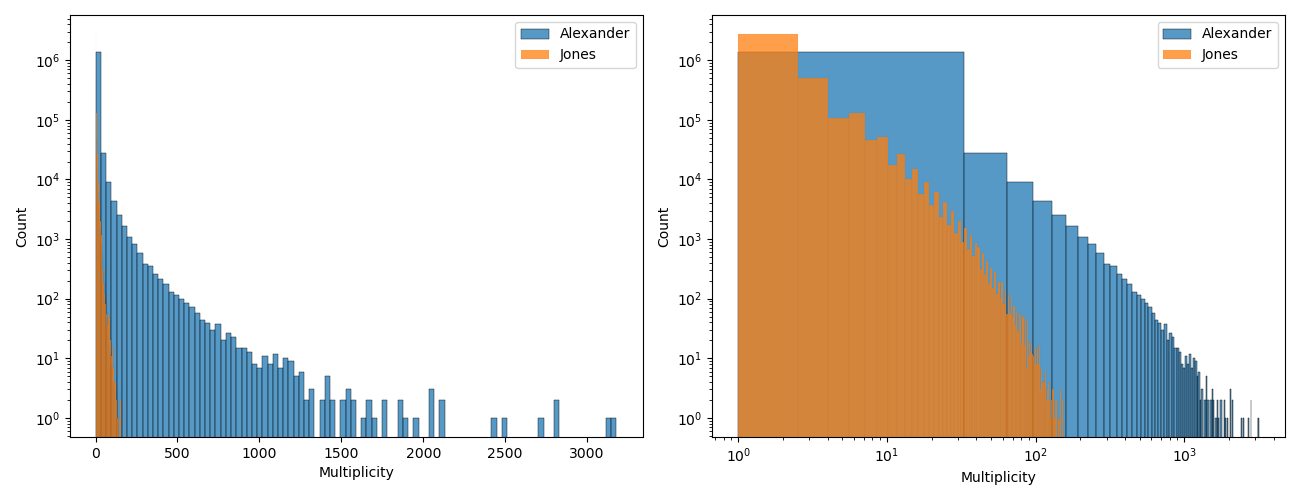}
	\caption{Illustration of a discriminative power of the Alexander and Jones polynomials of knots up to 17 crossings where x-axis represent the multiplicity and y-axis the count with  Alexander (blue) and Jones (orange).  }
	\label{fig:Mult}
\end{figure}

Since the it was shown that Khovanov homology distinguishes the unlink Hedden--Ni (2013), Batson--Seed (2015), trefoils Baldwin--Sivek (2018), yhe Hopf link  Baldwin--Sivek--Xie (2018), the connected sum of two Hopf links, the torus link $T(2,4)$ Xie--Zhang (2019),  split links Lipshitz--Sarkar  (2019), the torus link $T(2,6)$ Martin (2020), $L6n1$ Xie--Zhang (2020), $L7n1$, the connected sum of a trefoil and Hopf link  Li--Xie--Zhang (2020), and the figure eight $4_1$ \cite{baldwin2020khovanov}. Kanenoby knots \cite{kanenobu1986infinitely} fail to be distinguished by Khovanov and many other homology theories \cite{hedden2018geography}. In the case of the Jones polynomial, even detecting the unknot is still open for the Jones polynomial, although it has been checked for all knots up to 24 crossings \cite{tuzun2018verification, tuzun2021verification}. Note that there exists an infinite family of links with the trivial Jones polynomial \cite{eliahou2003infinite}. Although we know that there are infinite families of knots with the same Jones but distinct Alexander polynomial \cite{kanenobu1986infinitely} and while there are no infinite families of alternating knots with the same Alexander polynomial, they do exists for non-alternating knots \cite{kauffman2017infinitely}.

As an alternative to listing examples of knots where one invariant does better than the other and vice versa \cite{jablan2007linknot} we raise the \emph{botany} question \cite{hedden2018geography} in the context of data analysis.
\begin{center}
	\textsc{ Given a knot, how many distinct knots share its  invariant?}
\end{center}

Since simple search shows that among all knots with up to 17 crossings there are 2434 that have the same Alexander as the unknot, and  1668 the Alexander polynomial of a  trefoil we analyze the uniqueness and the size of non-unique classes of several polynomial invariants. In particular, the multiplicity of a polynomial refers to the number of knots up to certain crossing number sharing the same polynomial. If a polynomial has multiplicity $1$, it means that it is a polynomial of a single knot in a given collection. Such a value is called unique. On the other hand, polynomial has multiplicity $n>1$ if it is realized by $n$ knots from the considered collection. The count is the number of classes of knots with given multiplicity. For example there are 3 classes of 23 knots up to 15 crossings that have the same Alexander, Jones and HOMFLYPT polynomials. 

The Figure \ref{fig:Mult} compares  multiplicity (how many knots share the same polynomial) with the count (how many classes) indicated for Alexander by blue and for Jones by orange bars. Notice that on the left of Figure \ref{fig:Mult} the Alexander polynomial dominates Jones. The right panel shows the same data in case when the x-axis is on a logarithmic scale which provides finer details. Note that Jones polynomial have much more unique values compared to Alexander. Note that there are about 100 knots with Jones multiplicity around 100 but there are just under 100 Alexander polynomials with multiplicity more than 1000! 

The Table~\ref{tab:DistUni} provides more quantitative information, namely the percentage of unique (polynomials with multiplicity $1$) and distinct polynomials which include all unique polynomials and also those which are achieved by more than one knot for Alexander, Jones, HOMFLYPT, and Khovanov polynomials as well as all meaningful pairwise combinations: Alexander and Jones, Alexander and Khovanov, and HOMFLYPT and Khovanov. 

\begin{table}[h]
	\centering
	\begin{tabular}{|c|c|c|c|c|c|c|c|}\hline 
		& $\Delta$ & $V$   & $P$   & $Kh$  &  $\Delta \& V$  & $\Delta \& Kh$ & $P \& Kh$\\ 
		\hline\hline
		Unique $\%$     & 11.2     & 26.9  & 53.7  & 38.7   &  53.4 & 61.3 & 61.5 \\ \hline
		Distinct $\%$   & 24.6     & 47.6  & 72.0  & 59.4   &  71.6 & 77.4 & 77.6  \\ \hline
	\end{tabular}
	\caption{The percentage of knots up to 15 crossings 
		with distinct or unique value of the polynomial invariant or a pair of them for Alexander ($\Delta$), Jones ($V$), HOMFLYPT ($P$), and Khovanov ($Kh$) polynomials.  }
	\label{tab:DistUni}
\end{table}

\begin{figure}[ht!]
	\centering
	\includegraphics[width=0.85\linewidth]{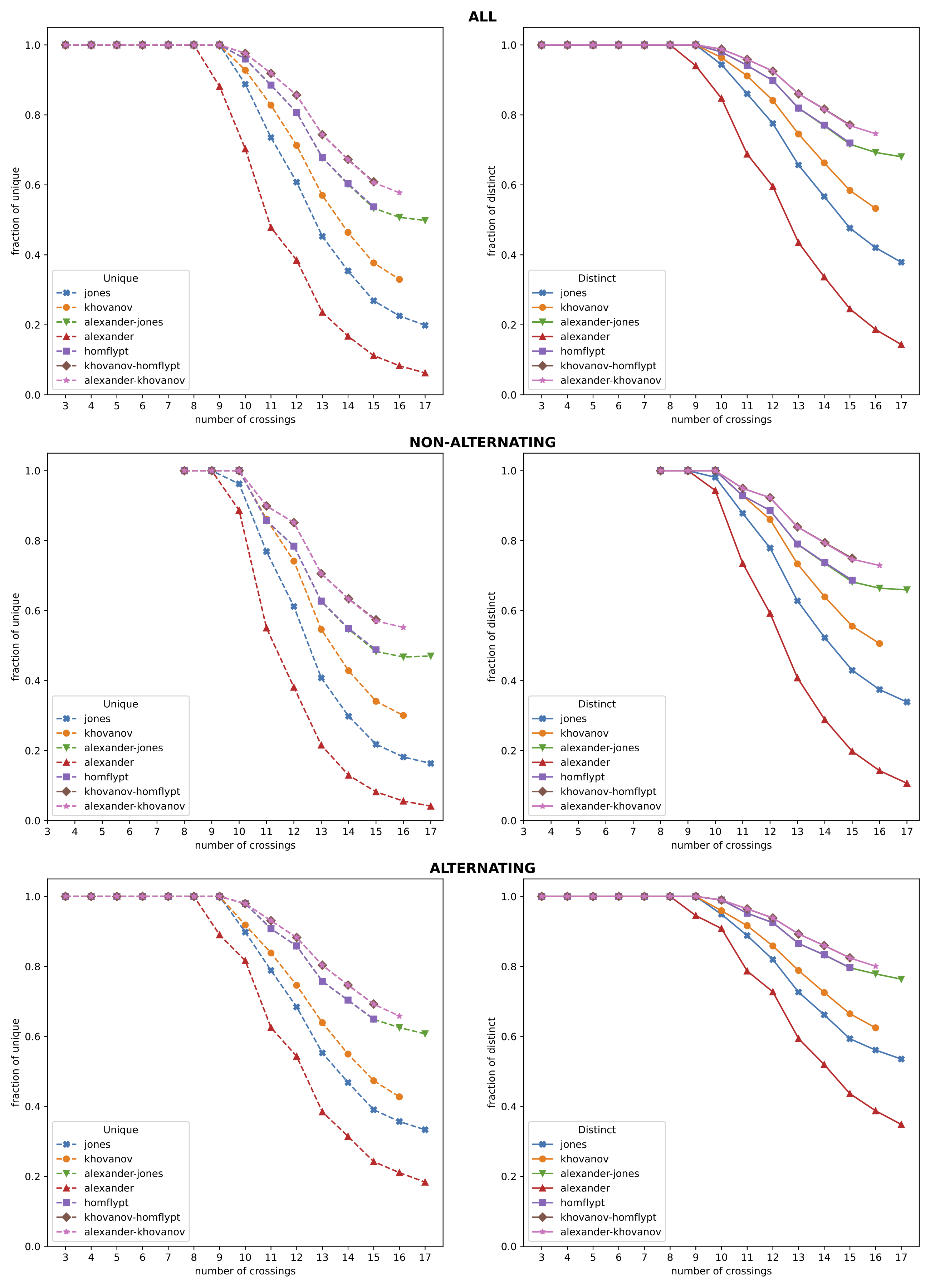}
	\caption{Graphs showing the percentage of unique (left) and distinct (right) values of a single or a pair of polynomial invariants across the filtration by the crossing number for all knots (top), just non-alternating (middle) or just alternating knots (bottom row).  }
	\label{fig:UniqD}
\end{figure}

It is important to point out that there is no inconsistency with Alexander polynomial detecting only 11.3\%, Jones just 26.9\% while together they achieve almost the best result of 53.4\% unique values. This is due to the multiplicity (the size of the class of knots with the same value of either Jones or Alexander) and says that they are, informally speaking, detecting different information so that combined, they are much better than alone. Actually this is a partial explanation why they perform almost as well as possible- where the upper bound is given by the discriminative power of the HOMFLYPT polynomial. 

Notice the \emph{less-than-random} performance of all polynomials except for the HOMFLYPT, and curious effectiveness of Alexander combined with Jones - better than Khovanov and almost as good as HOMFLYPT, at least for knots up to 15 crossings- observation that deserves theoretical considerations. 
Though we conjecture that the ratio of the number of unique polynomials to the number of distinct knots converges to zero as the number of crossings increases for all invariants that do not distinguish mutants\footnote{Pointed out independently by Danniel Tubbenhauer (University of Sydney) and Jason Cantarella (UGA).}, we also expect that the probability of an invariant distinguishing two random knots approaches one\footnote{Noticed by Ernesto Lupercio (CINESTAV).}. Figure \ref{fig:UniqD} illustrates the decline of the percentage of knots with bounded number of crossings that can be distinguished by the invariant (unique, left) and the percentage of distinct invariants with respect to the total number of different knots (distinct, right) behaves as the number of crossings increases. The information in Figure \ref{fig:UniqD} extends the reach of data in Table \ref{tab:DistUni} to show the declining trend but also to emphasize lower discriminative power on non-alternating than alternating knots, Figure \ref{fig:UniqD} (middle, bottom), that yields the trend for all knots Figure \ref{fig:UniqD}
(top).  
Our random knots data set of Jones polynomials contained about 25\% of distinct values.

\section{Fox conjecture: illustrated and extended }\label{ssec:Fox}

The Fox trapezoidal conjecture \cite{fox1962some}, which characterizes the behavior of the coefficients of the Alexander polynomials, has remained open since the 1960s.

\begin{conj}[Fox trapezoidal conjecture \cite{fox1962some}]
	Let $K$ be an alternating knot with the Alexander polynomial given by 
	$\displaystyle \Delta_K(t) = \sum_{j=0}^{2n}(-1)^ja_jt^{2n-j}$
	for $a_j>0.$ Then   $$a_0 < a_1 < \ldots %< a_{n-m-1} 
	< a_{n-m} = \ldots  = a_{n+m} > a_{n+m+1} > \ldots > a_{2n}.$$
\end{conj}

Additionally, Hirasawa and Murasugi \cite{hirasawa2014stability} conjectured that $m \leq|\sigma(K)|/2$ implying that the length of the stable part of the Alexander polynomial is bounded above by signature. 
\begin{figure}[h]
	\centering
	\begin{subfigure}[b]{0.45\textwidth}
		\centering
		\includegraphics[width=0.8\textwidth]{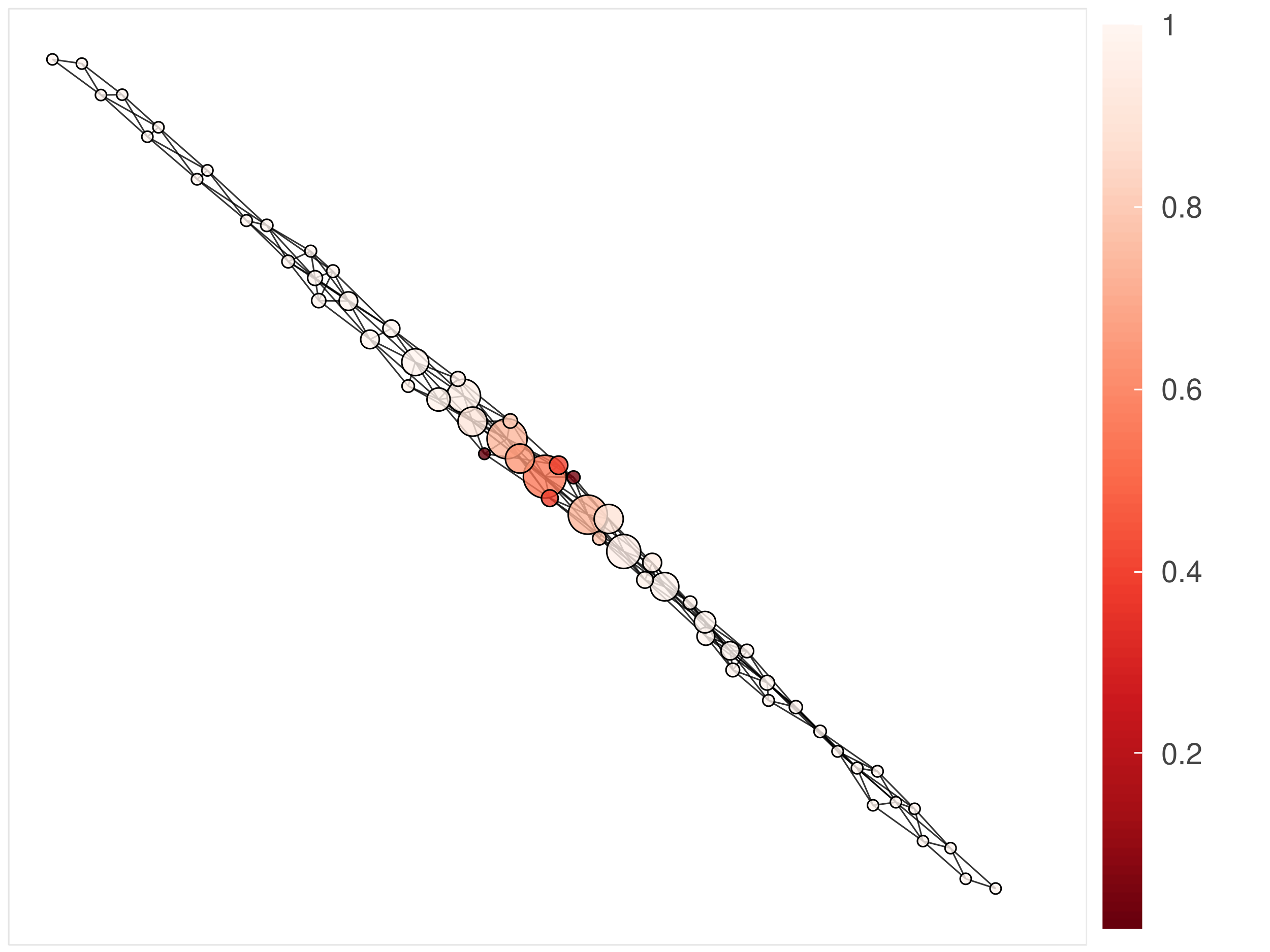}
		\caption{}\label{fig:alexlin}
		\hfill
	\end{subfigure}
	\begin{subfigure}[b]{0.45\textwidth}
		\centering
		\includegraphics[width=0.8\textwidth]{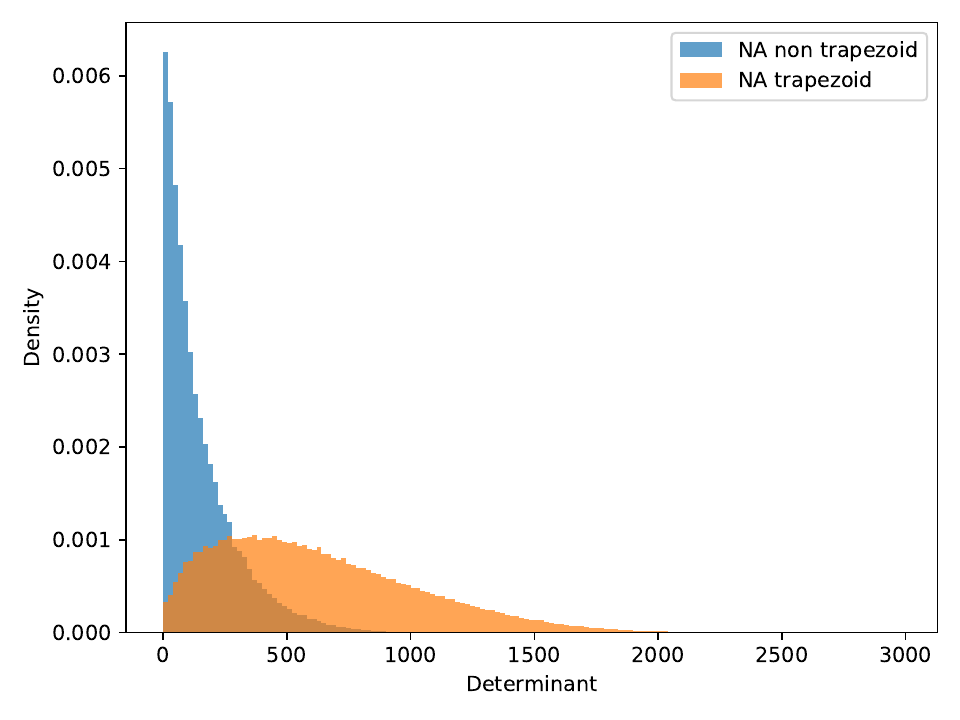}
		\caption{}
	\end{subfigure}
	\caption{(A) Ball Mapper of the Alexander polynomial for all all non-alternating  knots up to 17 crossings colored by the fraction of knots for which the Fox conjecture holds. (B) The distribution of the determinant values for the non-alternating knots up to 17 crossings on which the Fox conjecture holds (blue) and does not hold (orange).}
	\label{fig:FOX}
\end{figure}

The conjecture is known to be true for many classes of alternating knots: alternating pretzel knots \cite{parris1978pretzel},  two-bridge knots \cite{hartley1979two},  alternating algebraic knots \cite{murasugi1985alexander},  genus two alternating
knots \cite{ozsvath2003absolutely}.
In 2024, Fox conjecture was proven for special alternating knots \cite{hafner2024alexander, hafner2024log}, their plumbings and additional classes \cite{azarpendar2024fox}.

A closer look at the Alexander polynomial data for knots up to 17 crossings shows that the Fox conjecture holds for 6,346,032 out of 7,494,022 non-alternating knots, which is more than $85\%$. However, only 233,192 are trapezoids and not just triangles, which is just over $3\%$.  1829696 NA knots up to 17 crossings with signature 0.

Ball Mapper for the Alexander polynomial of all non-alternating knots with up to 17 crossings colored by the fraction of knots for which it holds, Figure \ref{fig:FOX} shows that such non-alternating knots are concentrated in central nodes. In \cite{dlotko2024mapper} we show that the values of the determinant are lowest in central nodes but increase in absolute value as we move away from the center. Our hypothesis, illustrated in the plot in Figure \ref{fig:FOX}, is that non-alternating knots for which Fox conjecture does not hold (orange) have smaller values of the determinant.

%%%%%%%%%%%%%%%%%%%%%%%%%%%%%
%\section{Ball Mapper}\label{sec:BM}
\section{Ball Mapper applications in knot theory}\label{sec:BM}

\begin{figure}[ht]
	\centering
	\begin{subfigure}{0.32\textwidth}
		\includegraphics[width=\textwidth]{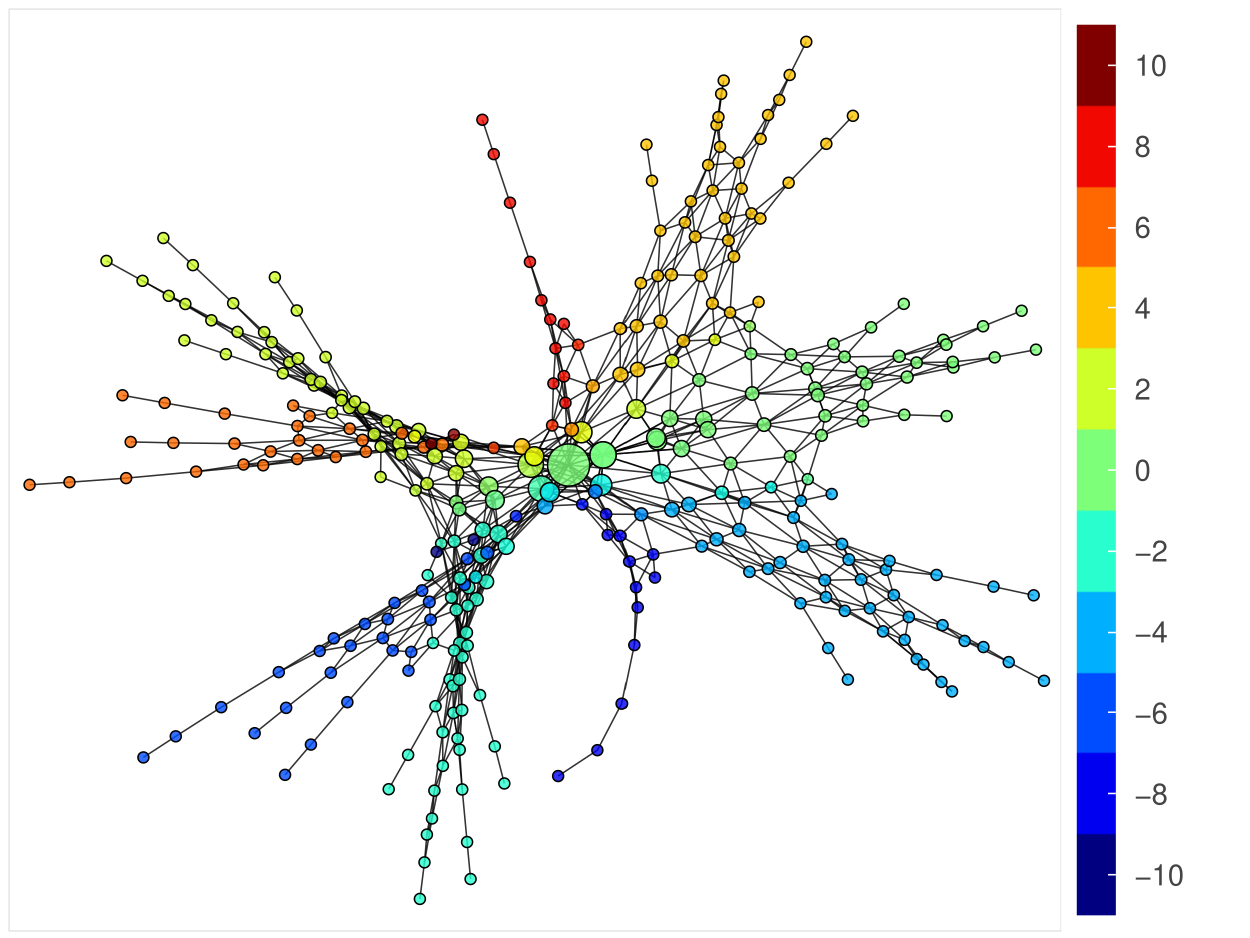}
		\caption{Jones data.}
		\label{fig:secondJ}
	\end{subfigure}
	%\hfill
	\begin{subfigure}{0.32\textwidth}
		\includegraphics[width=\textwidth]{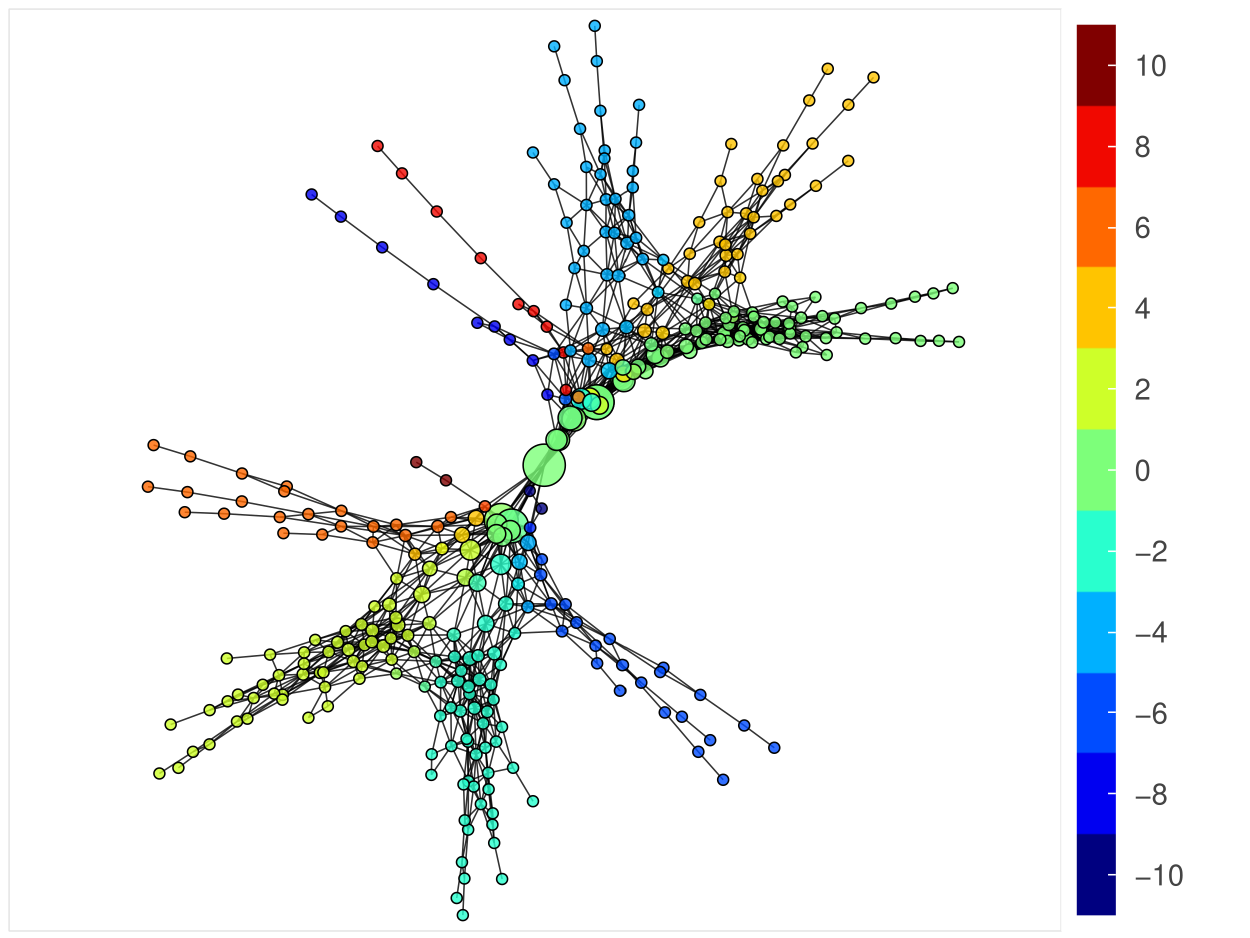}
		\caption{Alexander and Jones data.}
		\label{fig:thirdAJ}
	\end{subfigure}
	%\hfill
	\begin{subfigure}{0.32\textwidth}
		\includegraphics[width=\textwidth]{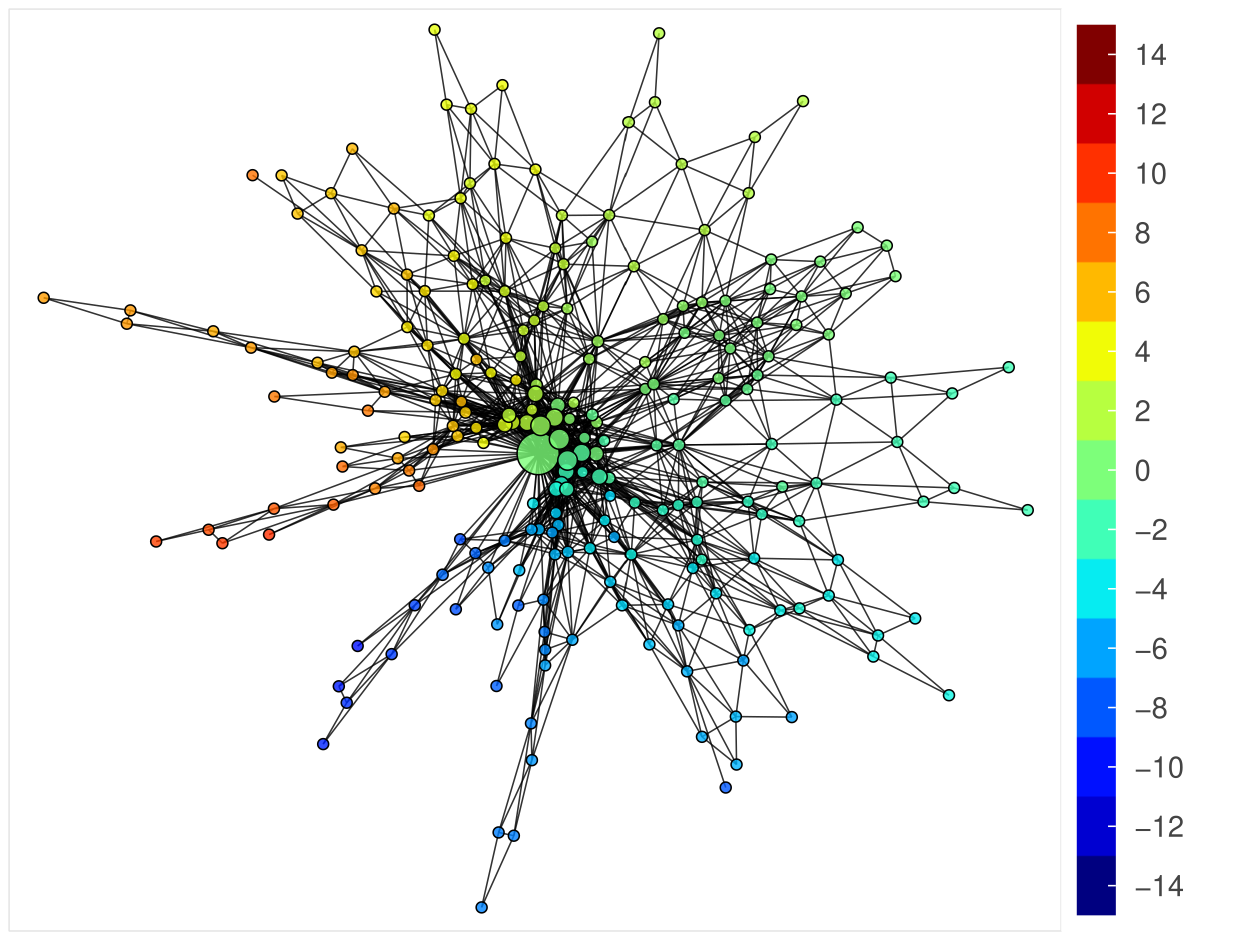}
		\caption{HOMFLY-PT data.}
		\label{fig:thirdHPT}
	\end{subfigure}
	
	\caption{Ball Mapper plots for polynomial data of all knots up to 15 crossings and their mirrors for the 
		%Alexander (A),
		Jones (A), HOMFLY-PT (C) and the combination of the Alexander and Jones data together (B) colored by signature.}
	\label{fig:BM_upto_15}
\end{figure}

The empirical confirmation of the stability of Ball Mapper graphs on Jones polynomial data (see~\cite{dlotko2024mapper}), and Alexander, HOMFLYPT and Khoavnov polynomials satisfy the same stability with respect to the crossing number filtration and the choice of parameter $\epsilon$. In this section we analyze properties and the star-like structure of BM graphs.  One can observe the dense collections of nodes in the center, the majority of the knots lies in the center of the Ball Mapper graph, 
% Figure \ref{fig:firstA},  
and there is a number of flares emanating from it. Note that the flares correspond nicely with the average value of the signature of knots in each node, Figure \ref{fig:BM_upto_15}.
%\ifdefined\ShowPictures
\begin{figure}[h!]
	\centering
	\begin{subfigure}[b]{0.35\textwidth}
		\includegraphics[width=0.9\textwidth]{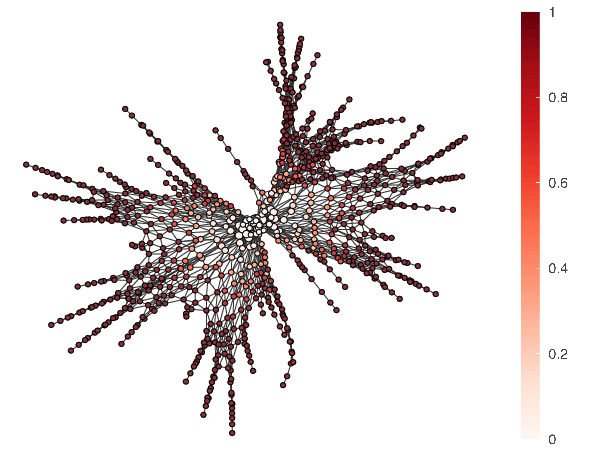}
		\caption{}
		\hfill
	\end{subfigure}
	\begin{subfigure}[b]{0.35\textwidth}
		\includegraphics[width=0.9\textwidth]{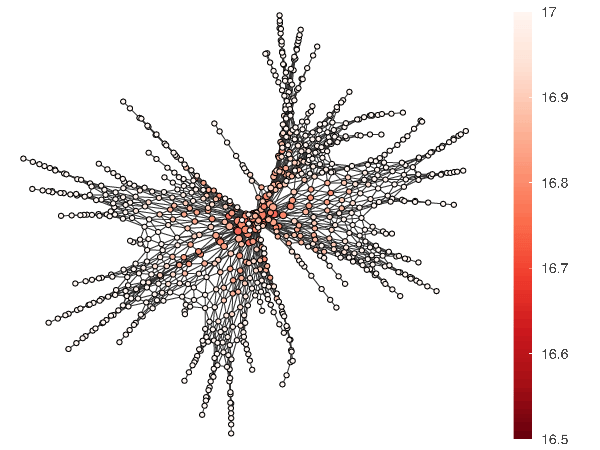}
		\caption{}
	\end{subfigure}
	\caption{Jones space of knots of up to 17 crossings colored by the fraction of alternating knots in each cluster (A),  and by the crossing number (B).}
	\label{fig:jones_upto_17Three}
\end{figure}
%\fi

Additionally, Figure \ref{fig:jones_upto_17Three} provides evidence that nodes in the center of the Ball Mapper  graph contain a mixture of alternating and non-alternating knots with a range of crossing numbers. On the other hand, flares consist mostly of the points corresponding to alternating knots with the highest crossing number in the data set. 

\begin{figure}[h!]
	\centering
	\begin{subfigure}[b]{0.3\textwidth}
		\includegraphics[width=\textwidth]{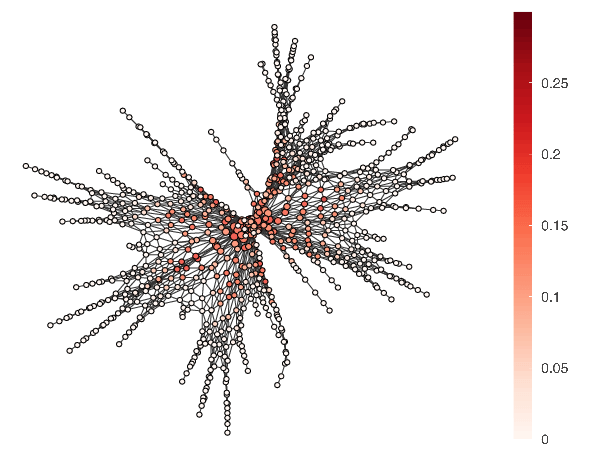}
		\caption{All 16 crossings.}
	\end{subfigure}
	\begin{subfigure}[b]{0.3\textwidth}
		\includegraphics[width=\textwidth]{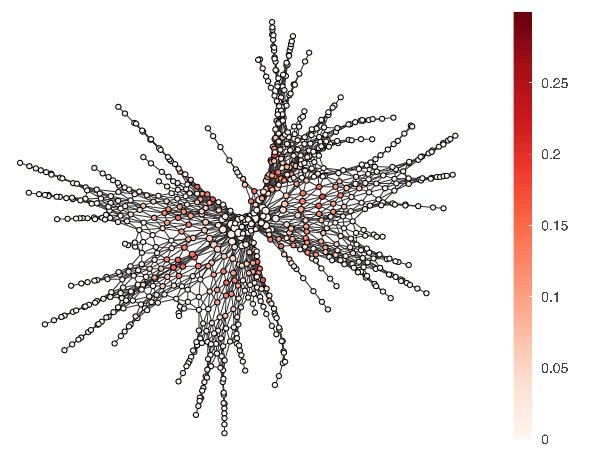}
		\caption{16 crossings alternating}
	\end{subfigure}
	\begin{subfigure}[b]{0.3\textwidth}
		\includegraphics[width=\textwidth]{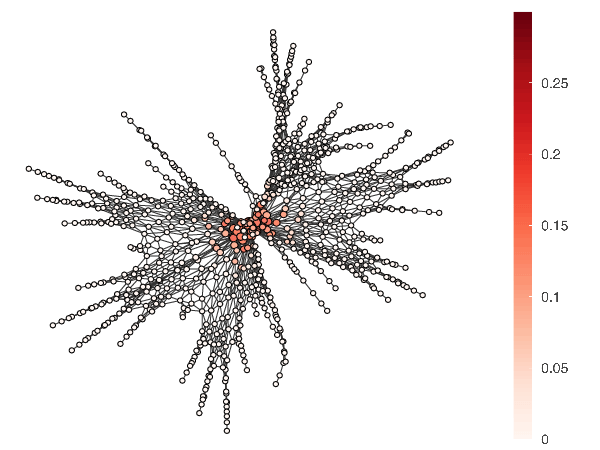}
		\caption{16 crossings non-alternating}
	\end{subfigure}
	\caption{ Ball Mapper graph for the Jones polynomial data of knots up to $17$ crossing colored by the fraction of knots in each cluster with $16$ crossings. }
	\label{fig:16in17}
\end{figure}

Following up on this insight, Figure \ref{fig:16in17}, shows the distribution of points corresponding to $16$ crossing knots within the BM of Jones data for $17$ crossing knots. In summary, 16 crossing alternating knots are underrepresented in the center and very ends of the flares while non-alternating $16$ crossing knots belong to central nodes, and similar structures appear when visualizing data corresponding to crossing number equal to $m$ within data of all knots with up to $n>m$ crossings for $n 
\leq 17.$

\begin{figure}[h!]
	\centering
	\begin{subfigure}[b]{0.4\textwidth}
		\centering
		\includegraphics[width=\textwidth]{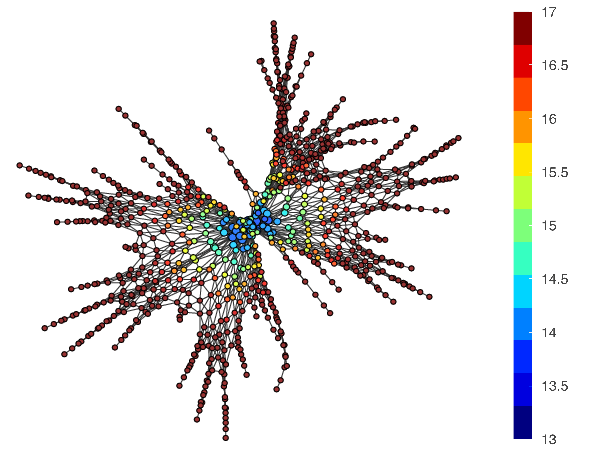}
		\caption{Jones span $\epsilon=100$}\label{spanJones}
	\end{subfigure}
	\hspace{1cm}
	\begin{subfigure}[b]{0.4\textwidth}
		\centering
		\includegraphics[width=\textwidth]{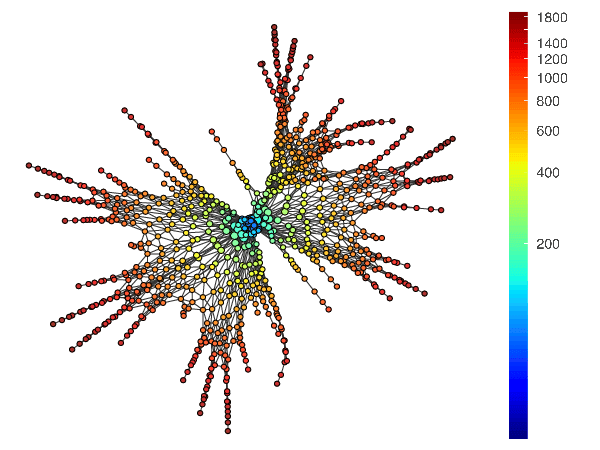}
		\caption{Determinant (logscale) $\epsilon=100$}\label{detJones}
	\end{subfigure}
	\caption{Ball Mapper of Jones polynomial data of knots up to $17$ crossings colored by the determinant (A) and the span of the Jones polynomial (B). }
	\label{fig:17span}
\end{figure}

\begin{figure}[ht]
	\centering
	% \begin{subfigure}{0.45\textwidth}
		%     \includegraphics[width=\textwidth]{MM_alex.png}
		%     \caption{Alexander data}
		% \end{subfigure}
	% \hfill
	\begin{subfigure}{0.45\textwidth}
		\includegraphics[width=\textwidth]{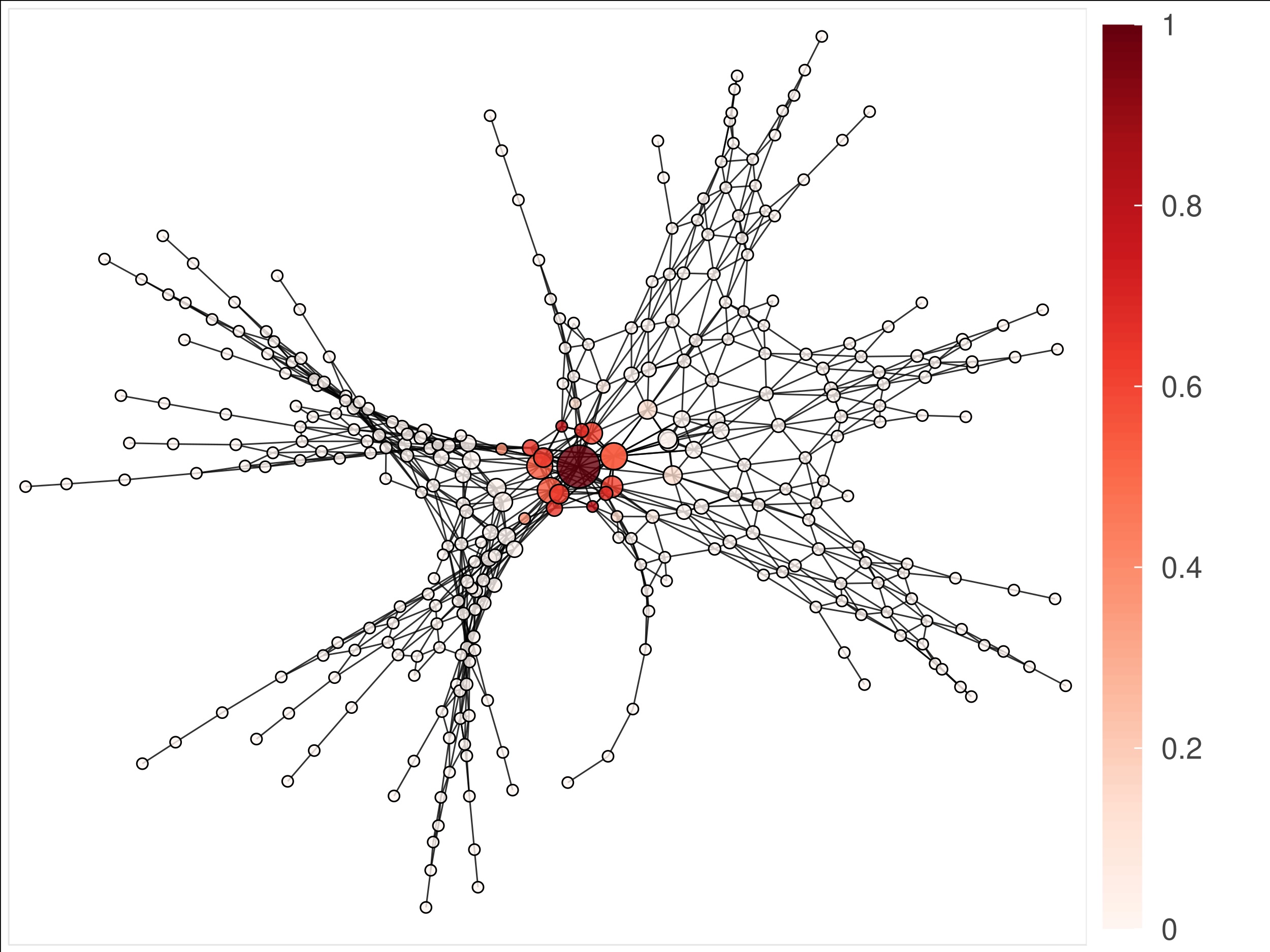}
		\caption{Jones data.}
		\label{fig:secondJJ}
	\end{subfigure}
	\hfill
	\begin{subfigure}{0.45\textwidth}
		\includegraphics[width=\textwidth]{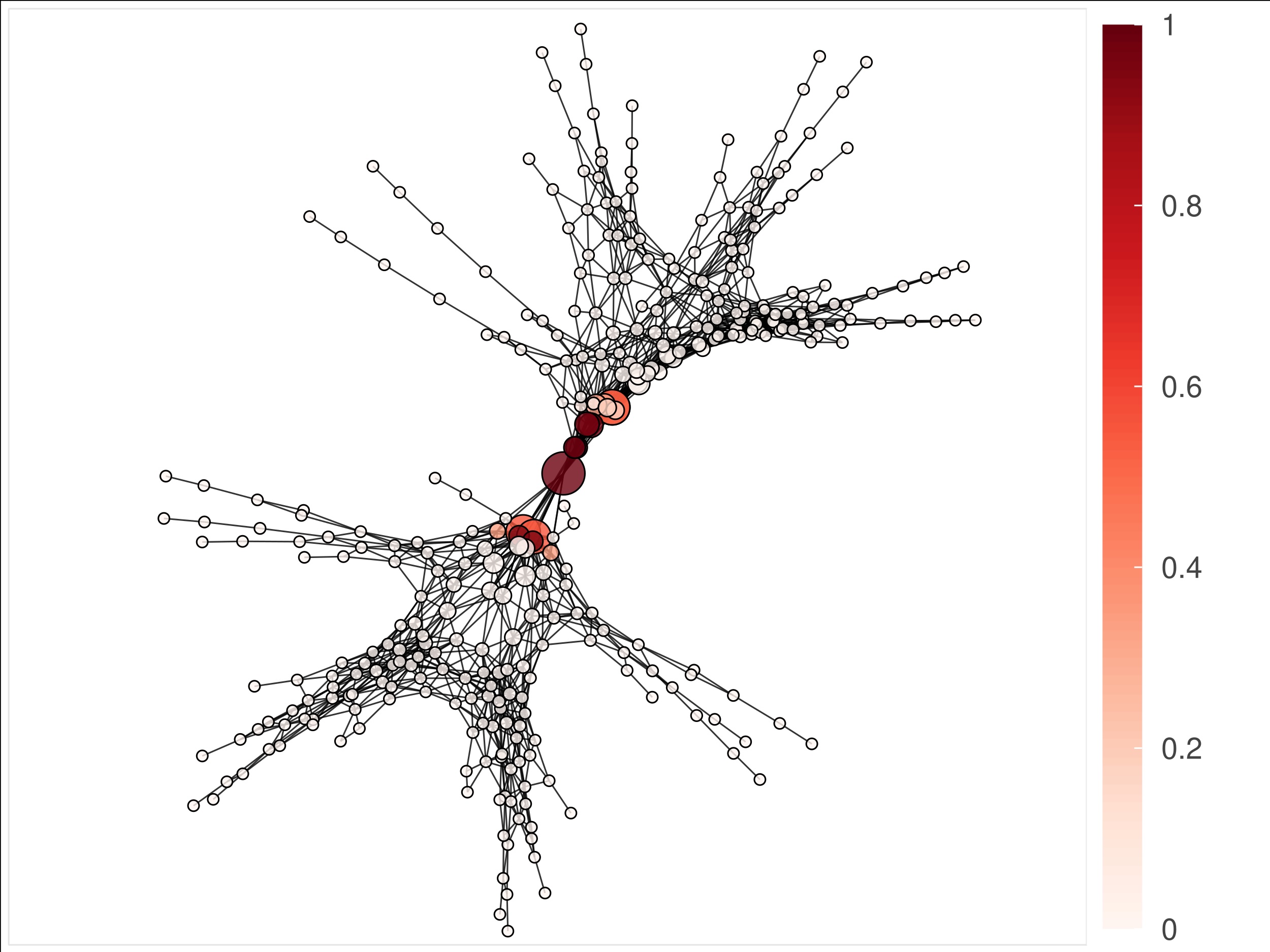}
		\caption{Alexander and Jones data.}
		\label{fig:thirdAJJ}
	\end{subfigure}
	\hfill
	\begin{subfigure}{0.45\textwidth}
		\includegraphics[width=\textwidth]{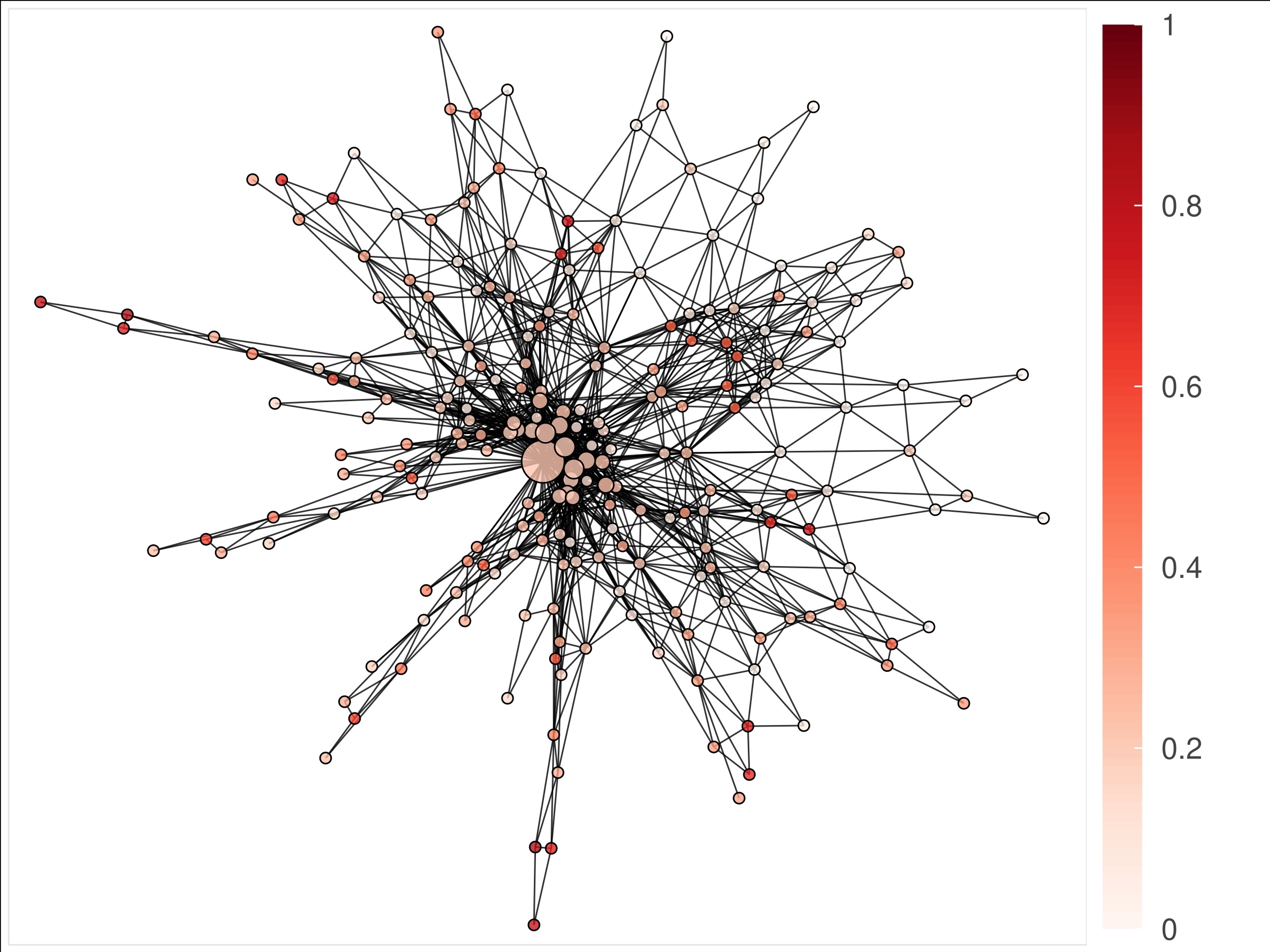}
		\caption{HOMFLY-PT data.}
	\end{subfigure}
	\hfill
	\begin{subfigure}{0.45\textwidth}
		\includegraphics[width=\textwidth]{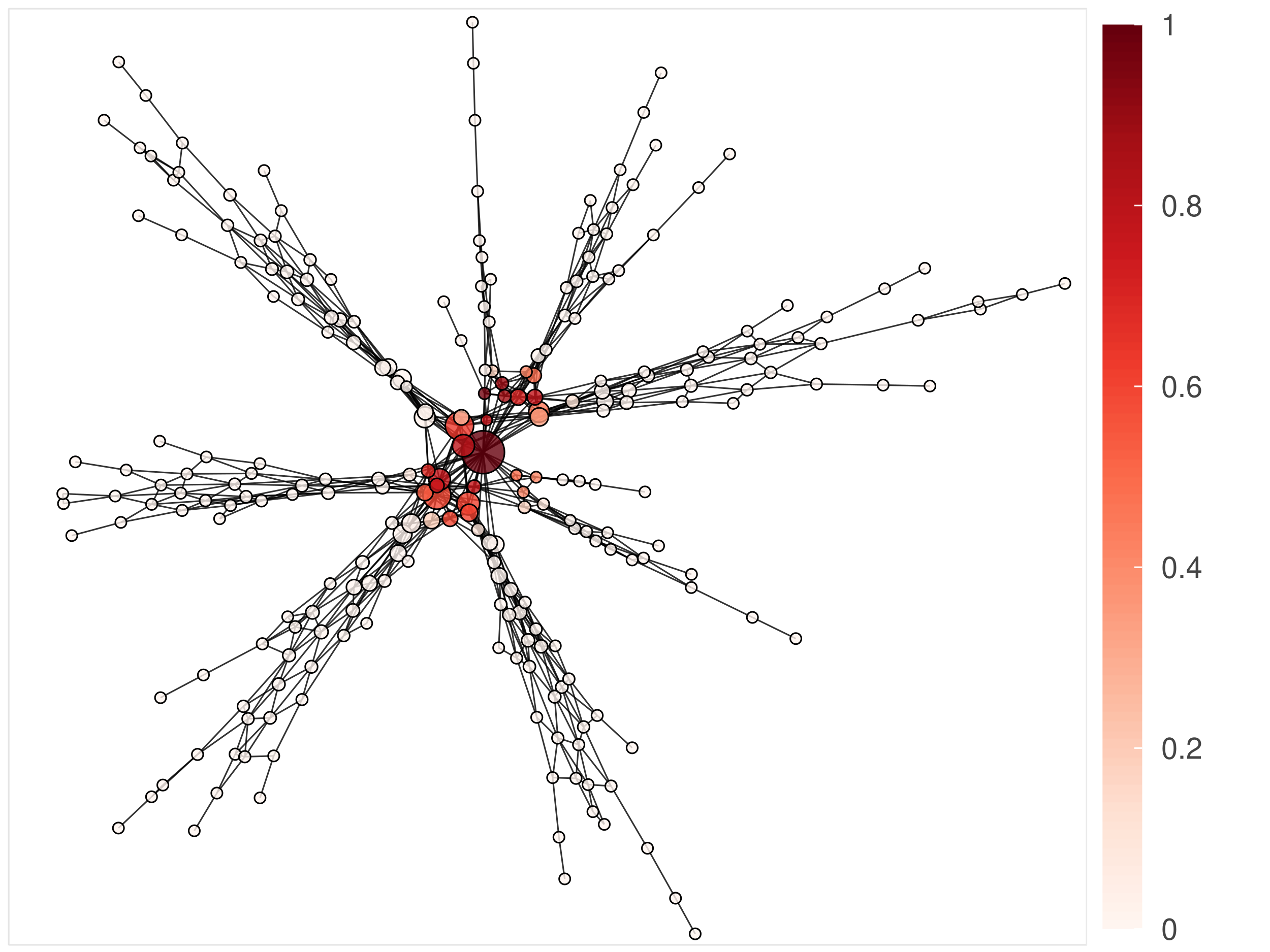}
		\caption{KH data.}
	\end{subfigure}
	
	\caption{Ball Mapper graphs for polynomial data of all knots up to 15 crossings, with mirrors for the %Alexander (A), 
		Jones (A), HOMFLYPT (C), Khovanov (D) polynomials and the combination of the Alexander and Jones data together (C). Colored by the fraction of knots in each node that are covered by the central node in Jones.}
	\label{fig:BM_upto_15_MappingMappers}
\end{figure}

Ball Mapper graphs \cite{datasets}  suggest that all of the insights from Jones data for all knots up to $17$ crossings carry over to Jones data of only alternating or just non-alternating knots, Figure \ref{fig:JonesNA}, with exactly $17$ crossings. Their Ball Mapper graphs remain star-like with the flares correlated with signature values. The central area in the Ball Mapper of alternating knots consist of a single node containing knots with their mirror images, hence, the average signature is zero. 
This implies that in the Ball Mapper of all knots the non-alternating knots are concentrated in central nodes that also have smaller span, see Figure \ref{spanJones}, and smaller determinant Figure \ref{detJones} implying the well known theorem that the span of the Jones polynomial of alternating knots determines their crossing number but provides just a strict lower bound for non-alternating knots. Note that the flares in the Ball Mapper are less prominent in just non-alternating data Figure \ref{fig:JonesNA}.

Ball Mapper graph for the Alexander polynomial data has linear structure and it is stable with respect to crossing number filtration and the choice of Ball Mapper parameter Figure \ref{fig:alexlin}.  
Ball Mapper graphs of HOMFLYPT polynomials of all knots up to $15$ crossings colored by the average signature of knots in each node are presented in Figure \ref{fig:thirdHPT}. The analogous image for Khovanov polynomial data for non-alternating knots up to $15$  crossings is shown on Figure \ref{fig:firstKH}.
Both datasets share the star-like structure with the Ball Mapper graph of the Jones polynomial and satisfy similar properties, in particular with respect to signature though flares of the HOMFLYPT Ball Mapper appear less differentiated and most knots belong to the highly connected center. However, Figure \ref{fig:BM_upto_15_MappingMappers} highlights a crucial difference between Ball Mapper graphs of the Jones, Alexander and Jones, and Khovanov data with respect to the HOMFLYPT polynomial. MappingMappers algorithm illustrated in this figure shows that the central node in the Ball Mapper of the Jones data correspond to better differentiated central nodes in the Alexander and Jones combined, and Khovanov data while HOMFLYT has a very different distribution of those knots Figure\ref{fig:BM_upto_15_MappingMappers}(c).

Inspired by results in Section \ref{sec:intuition} we also provide Ball Mapper of data consisting of a both Alexander and Jones data concatenated of all knots up to 15 crossings Figure \ref{fig:thirdAJ}. Notice the more elongated structure and better separation of central vertices with respect to either Alexander or Jones Ball Mappers. Yet, HOMFLYPT and the Alexander-Jones pair detect approximately the same percentage of knots Table \ref{tab:DistUni} and have similar decay curves with respect to the crossing number Figure \ref{fig:Mult}, indicating that the Ball Mapper structure is not always indicative of the discriminative power of the invariant since it clusters knots with the similar values of the invariant. Interactive Ball Mapper graphs and visualizations of KnotInfo data \cite{knotinfo} are available at \cite{dioscuri}.

%%%%%%%%%%%%%%%%%%%%%%%%%%%%%
%\section{Ball Mapper vs persistent PCA}\label{sec:BMvsPCA}
%\input{BMvsPCA}
\section{Ball Mapper vs Persistent PCA}\label{sec:BMvsPCA}

The first paper to integrate topological data analysis with standard data analysis techniques in knot theory introduces the concept of Persistent PCA and investigates the dimensionality of the associated point cloud \cite{levitt2019big}.  

\begin{figure}[h]
	% \centering
	% \begin{subfigure}[]{0.45\textwidth}
		\includegraphics[width=6cm]{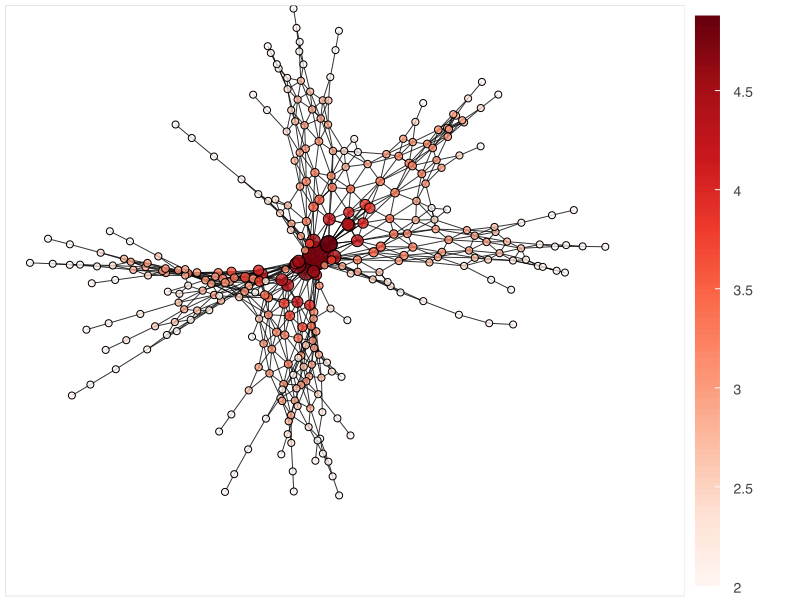} \hspace{1cm}
		\includegraphics[width=6cm]{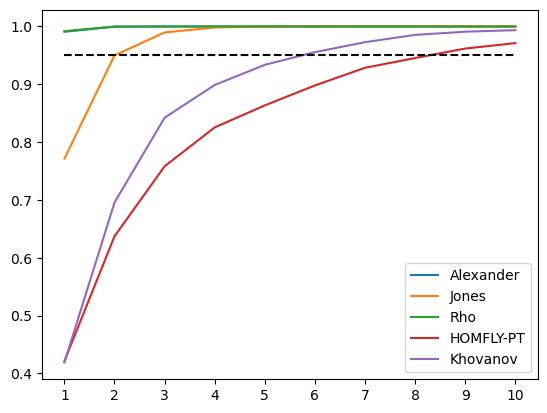}
		%      \caption{}
		% \end{subfigure}
	%     \begin{subfigure}[]{0.45\textwidth}
		% \input{BMvsPCAPCADimensionality.jpeg} 
		%        \caption{}\labe{fig:PCADimen}
		%     \end{subfigure}
	\caption{Ball Mapper graph of the Jones polynomial data for knots up to 17 crossings colored by the local dimension i.e. the dimension of each cluster (left). The cumulative explained variance plotted against the number of the principal components for the Alexander, Jones, $\rho$, HOMFLYPT and Khovanov polynomials of knots up to 15 crossings (right). }
	\label{fig:BMLocalDim}
\end{figure}

The dimensionality of the point cloud is defined to be the smallest value 
$\mathbf{d}$
for which the normalized explained variance of the first 
$\mathbf{d}$
principal components exceeds 95\%  across all considered crossing numbers.
In particular, they observe that the Alexander polynomial data is roughly 1-dimensional, while the Jones data is at least 3-dimensional. Our computations up to 15 crossings confirm results in \cite{levitt2019big} and we extend them to show that $\rho$ is 1-dimensional, Khovanov 6-dimensional and HOMFLYPT 8-dimensional, see Figure \ref{fig:BMLocalDim}b.  Jones Ball Mapper graph colored by the local dimension of each cluster Figure \ref{fig:BMLocalDim}a supports this claim. The majority of clusters is 3-dimensional with only the huge central nodes of local dimension five, which is still only about 10\% of the ambient dimension equal to 51 for 17 crossing knots.  

MappingMappers was used to visualize  relation between Alexander and Jones data \cite{dlotko2024mapper} and further used to rediscover the folklore theorem \cite{murasugi_knot_2008} stating that the Alexander and Jones polynomial determine the value of signature modulo four. In \cite{dlotko2024mapper} we show that the linear gradient on the Alexander data corresponds to the values of the determinant (zero in the center and increasing/decreasing along the flares) suggesting that the determinant might be related with the 1st persistent principal direction. 

\begin{table}[h]
	\centering
	\begin{tabular}{|c|c|c|c|c|c|c|}\hline 
		& $\Delta$ & $V$ & $P$ & $Kh$ &  $\rho$ 
		\\ \hline
		Pearson coeff.  & 0.9979 & -0.98 & 0.016  & -0.891& -0.02\\ \hline
	\end{tabular}
	\caption{The Pearson correlation coefficient between the value of the determinant and the first principal component of the polynomial data for Alexander, Jones, HOMFLYPT, Khovanov and $\rho$ polynomials of knots up to 15 crossings. }
	\label{tab:PCACorr}
\end{table}

Pearson correlation coefficients between the determinant  and the first principal component are included in Table \ref{tab:PCACorr} showing strong correlation for the Alexander and Jones, moderate for Khovanov and  essentially no correlation for HOMFLYPT and $\rho$.  The moderate correlation of the fist PCA of Khovanov's data stems from the correlation with knots of signature equal to zero and two which are the majority of the knots in our data and we do not expect it to persist.

\begin{figure}[h]
	\centering
	\includegraphics[width=0.95\linewidth]{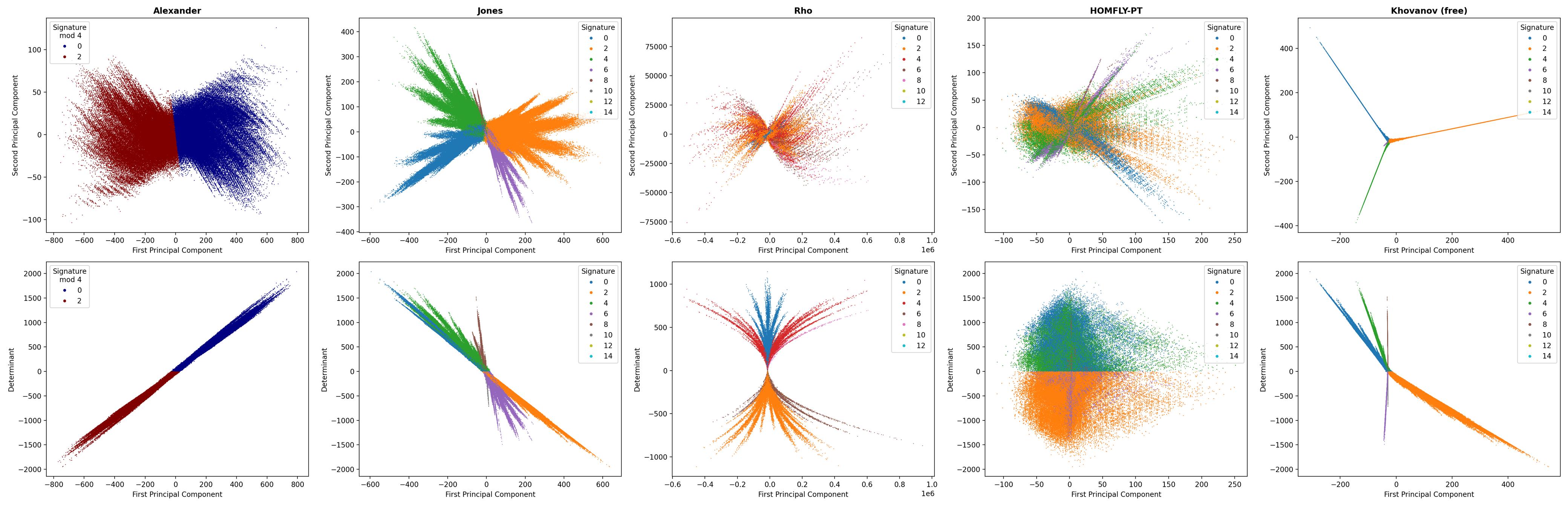}
	\caption{Top row: A PCA projection into 2 dimensions and bottom row: correlation between the two most significant principal components and the determinant for the Alexander, Jones, $\rho$, HOMFLYPT and Khovanov data for up to 15 crossing knots, from left to right. The Alexander graphs are colored by signature modulo 4 and the others by the value of signature. }
	\label{fig:DetPCA}
\end{figure}

Alexander Ball Mapper  data and the 2-dimensional projection indicate that the first principal component is essentially the determinant, and this is also true for the Jones polynomial.

%%%%%%%%%%%%%%%%%%%%%%%%%%%%%
%\section{sig}\label{sec:ssignature}
%\input{S-signature}
%%%%%%%%%%%%%%%%%%%%%%%%%%%%%
\section{Khovanov homology, s-invariant, and signature}\label{sec:ssig}
\label{sec:s_vs_sgnature}

In order to gain more insight into curious relation between the shape of Ball Mapper graphs for the Jones and  Khovanov data with signature witnessed by the monochormatic flares, we analyze the relation between the flares of these Ball Mappers and support of Khovanov homology. Khovanov homology is a bigraded homology theory which is supported in a number of diagonals of slope 2 (called the width of Khovanov homology. Knots whose Khovanov homology is supported on only two diagonals are called Khovanov thin. For example, alternating and quasialternating knots are homologically thin \cite{manolescu2007khovanov} but the opposite is not true \cite{greene2009homologically}. Moreover, it is known that signature and Rasmussen s-invariant coincide for alternating knots and that both provide a lower bound on the unknotting number, while s-invariant also bound the slice genus from below \cite{rasmussen2010khovanov}.

\begin{figure}[ht]
	\centering
	\begin{subfigure}{0.23\textwidth}
		\includegraphics[width=0.95\textwidth]{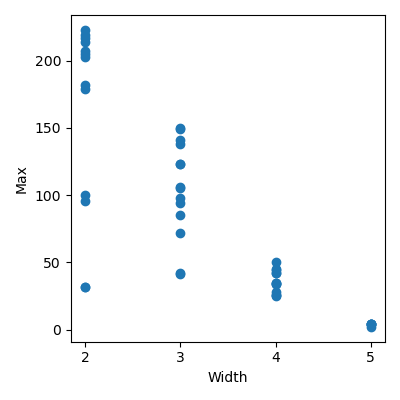}
		\caption{}\label{f:diagmax}
	\end{subfigure}
	\hfill
	\begin{subfigure}{0.23\textwidth}
		\includegraphics[width=0.95\textwidth]{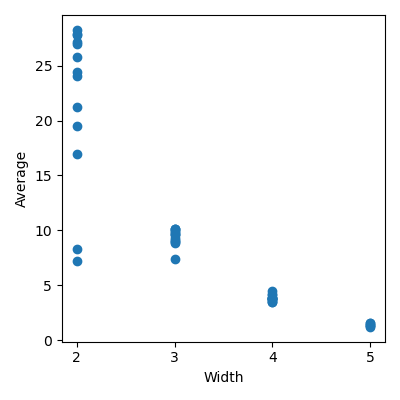}
		\caption{}\label{f:diagmean}
	\end{subfigure}
	\hfill
	\begin{subfigure}{0.23\textwidth}
		\includegraphics[width=0.95\textwidth]{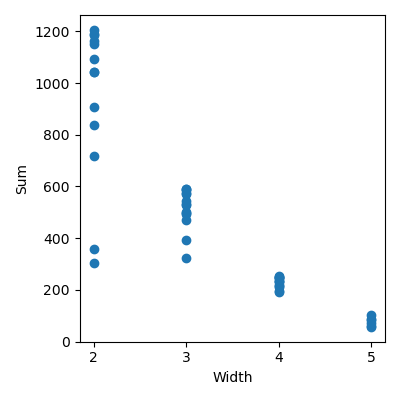}
		\caption{}\label{f:diagsum}
	\end{subfigure}
	\hfill
	\begin{subfigure}{0.23\textwidth}
		\includegraphics[width=0.95\textwidth]{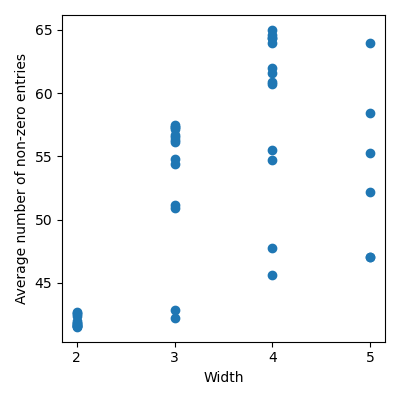}
		\caption{}\label{f:diag0}
	\end{subfigure}
	
	\caption{Statistics for the coefficients of Khovanov polynomials of all non-alternating knots up to 17 crossings grouped by diagonals: average maximal coefficient (A), average value of the coefficients (B), average sum of the coefficients (C), average number of non-trivial coefficients (D).}
	\label{fig:diags_statistics}
\end{figure}

Through the rest of this section we consider mainly Khovanov homology of 7.494.008 non-alternating knots up to 17 crossings, with about 3 million of them non-thin (they have more than 2 diagonals). 
Table \ref{tab:KhWidth} provides distribution of non-alternating knots with different width of Khovanov homology per crossing number, showing that thin non-alternating knots dominate the distribution though less so as the crossing number increases since the opportunities for higher width appear. 
Looking closely into simple statistics such as  maximal Figure \ref{f:diagmax}, average coefficients Figure \ref{f:diagmean}, and   the number of non-trivial entries Figure \ref{f:diag0} we notice that thin non-alternating knots have, on average, higher maximal and average, but fewer non-trivial coefficients. 

\begin{figure}[ht]
	\centering
	\begin{subfigure}{0.45\textwidth}
		\includegraphics[width=\textwidth]{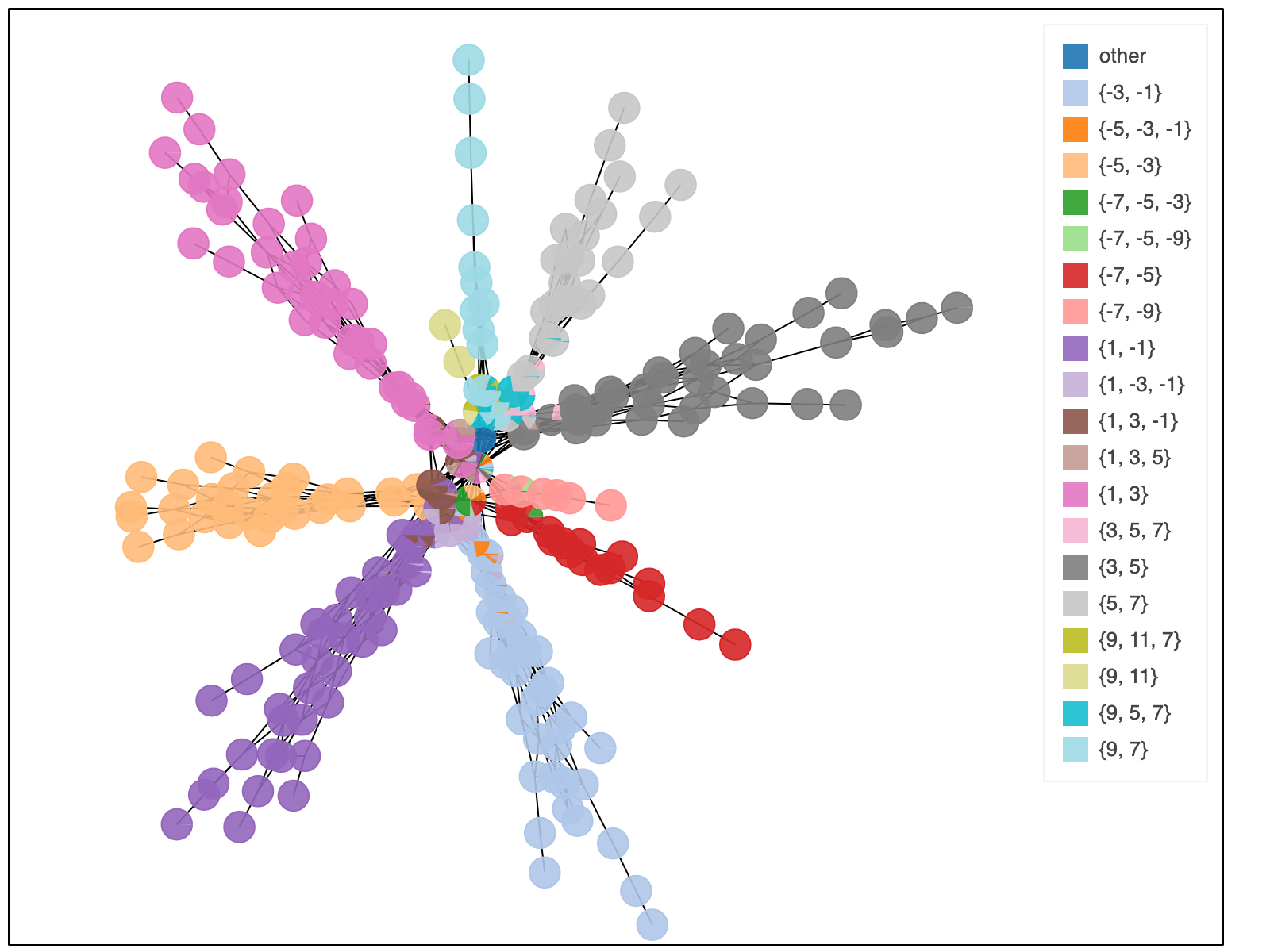}
		\caption{Khovanov homology data} \label{fig:BM_upto_15_diagKh}
		% \label{fig:firstA}
	\end{subfigure}
	\hfill
	\begin{subfigure}{0.45\textwidth}
		\includegraphics[width=\textwidth]{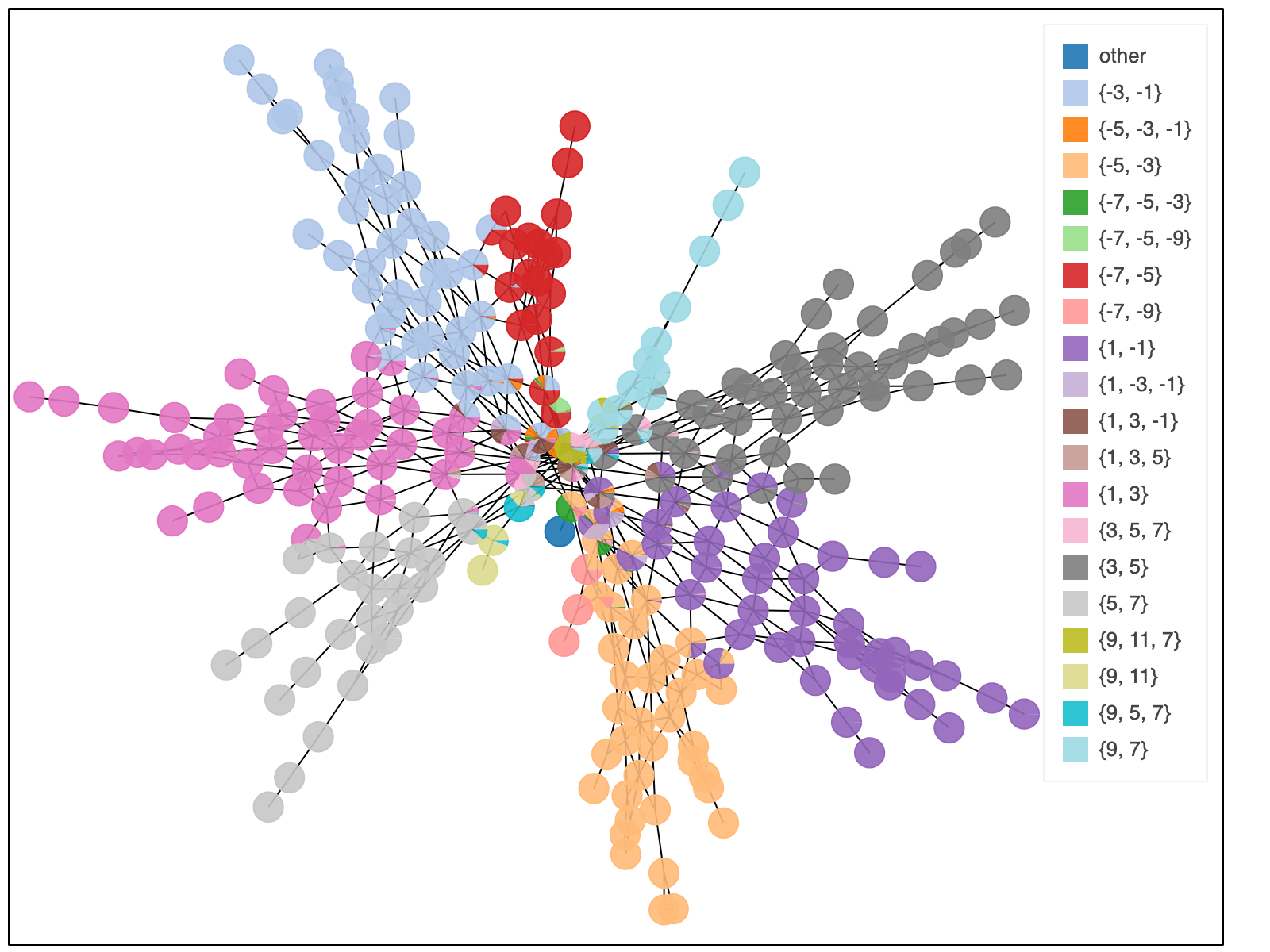}
		\caption{Jones polynomial data} \label{fig:BM_upto_15_diagJ}
	\end{subfigure}
	\caption{Khovanov homology and the Jones polynomial Ball Mapper graphs of knots upto 15 crossings with each node is represented by a piechart indicating the distribution of the different sets of diagonals in Khovanov homology in each node. }
	\label{fig:BM_upto_15_diag}
\end{figure}

Since the value of the s-invariant depends on non-trivial diagonals we test the hypothesis that the flares of  Khovanov Ball Mapper correspond to the diagonals Figure \ref{fig:BM_upto_15_diagKh}. Notice that all the vertices outside the center of  these ball mappers capture different combination of diagonals which correspond to nearby regions in our point cloud. Central vertices are colored by pie chats representing proportion of knots covered by that node in different diagonals. 

\begin{table}[h]
	\centering
	\begin{tabular}{|c|c|c|c|c|}  \hline
		\multirow{2}{*}{\parbox[c]{1.8cm}{\centering{Number of crossings}}} & \multicolumn{4}{c|}{Kh width} \\ \cline{2-5}
		& 2 & 3 & 4 & 5 \\ \hline
		% 8                                   & 2 & 1 & 0 & 0 \\ \hline
		% 9                                   & 9 & 2 & 0 & 0 \\ \hline
		% 10                                  & 41 & 12 & 0 & 0 \\ \hline
		% 11                                  & 185 & 53 & 0 & 0 \\ \hline
		% 12                                  & 854 & 272 & 0 & 0 \\ \hline
		13                                  & 4554 ($\sim$73\%)& 1678 ($\sim$27\%)& 4 & 0 \\ \hline
		14                                  & 23691 ($\sim$70\%)  & 9928 ($\sim$30\%)& 53 & 0 \\ \hline
		15                                  & 135054  ($\sim$67\%) & 66101($\sim$33\%) & 547 & 0 \\ \hline
		16                                  & 767345 ($\sim$63\%)& 437781 ($\sim$36\%)& 5480 & 0 \\ \hline
		17                                  & 4481503 ($\sim$60\%) & 2953589 ($\sim$40\%)& 58896  ($\sim$1\%)& 20 \\ \hline
	\end{tabular}
	
	\caption{The distribution of width of even Khovanov homology of all non-alternating knots from 13 to 17 crossings.}
	\label{tab:KhWidth}
\end{table}

This observation provides a natural and satisfying explanation of the shape of Khovanov Ball Mapper in terms of its structure/width and raises another question: why does the Jones polynomial Ball Mapper exhibit a similar coloring Figure~\ref{fig:BM_upto_15_diagJ} by different sets of diagonals. 

Recall that Khovanov homology \cite{bar2002khovanov, khovanov1999categorification} is a categorification of the Jones polynomial, which means that the $j$th coefficient of the Jones polynomial are equal to the alternating sum of Khovanov coefficients $Kh$ in $j$th quantum degree:
\begin{equation}\label{decat}
	% V(q)=Kh(-1, q) 
	% Kh (L)(q,-1)={\hat {J}}(L)(q):=(q+q^{-1})J(L)(q^{2})
	{Kh} (L)(q,-1)=(q+q^{-1})V(L)(q^{2}) .
\end{equation}

Although the decategorification \eqref{decat} results in a loss of information or data compression in data science terms, the structure of the diagonals in Khovanov homology is preserved and apparent in Jones Ball Mapper. Understanding the precise underlying mechanisms of this phenomenon would be interesting and aligned with questions posed in \cite{ davies2021advancing, craven2024illuminating}.

\begin{figure}[h]
	\centering
	\begin{subfigure}{0.32\textwidth}
		\includegraphics[width=\textwidth]{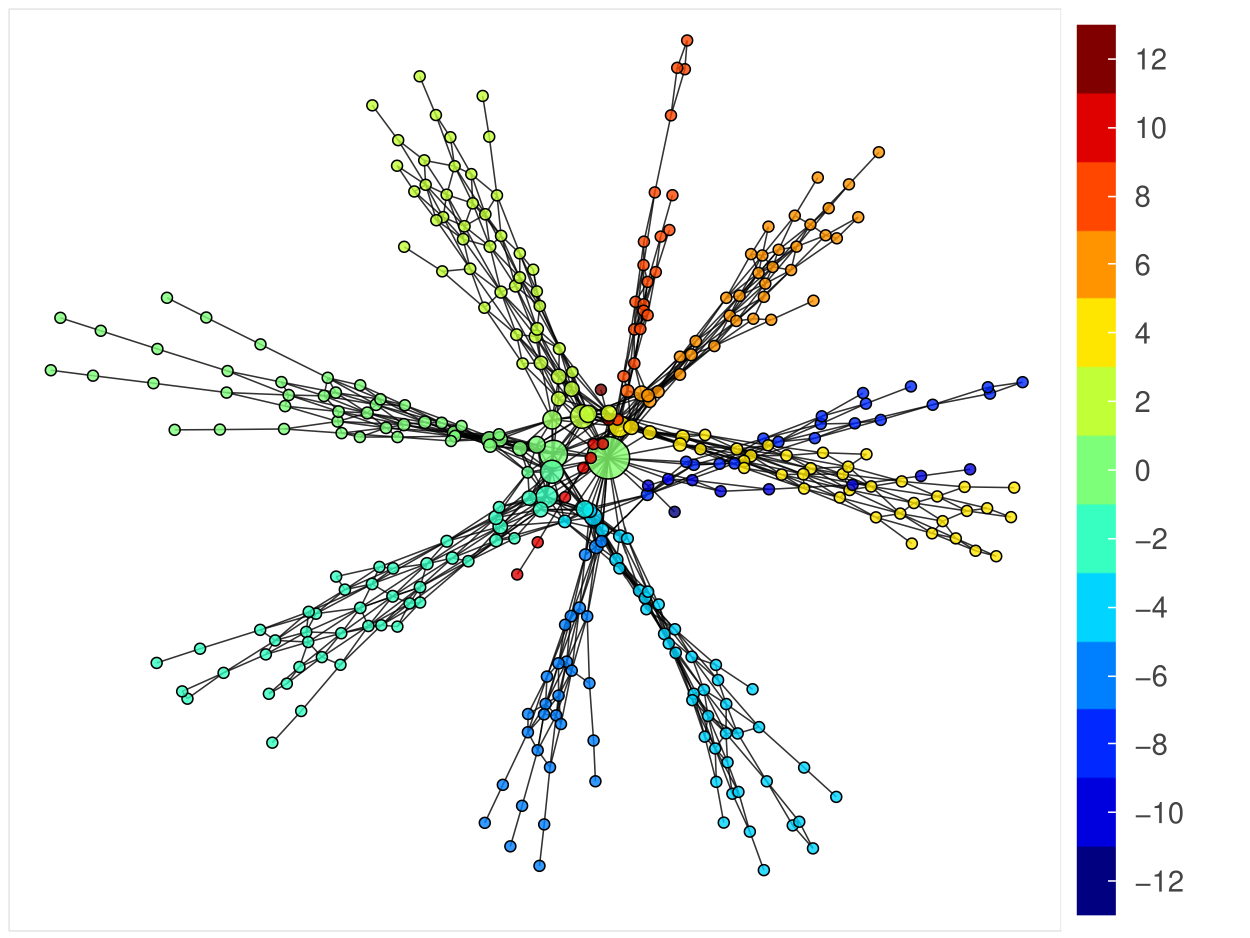}
		\caption{Khovanov all}
		\label{fig:firstKH}
	\end{subfigure}
	\hfill
	\begin{subfigure}{0.32\textwidth}
		\includegraphics[width=\textwidth]{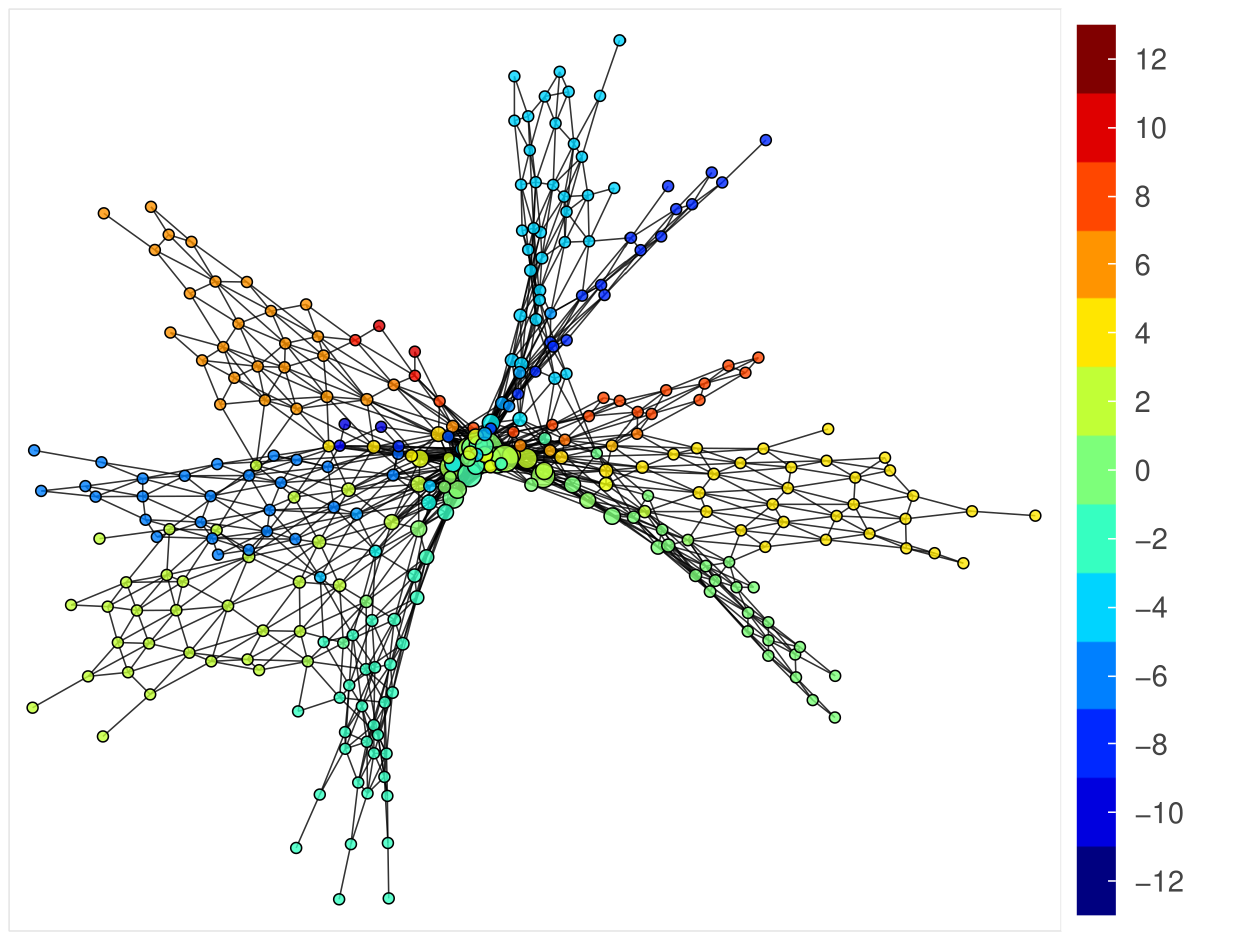}
		\caption{Jones for non-alternating}
		\label{fig:JonesNA2}
	\end{subfigure}
	\hfill
	\begin{subfigure}{0.32\textwidth}
		\includegraphics[width=\textwidth]{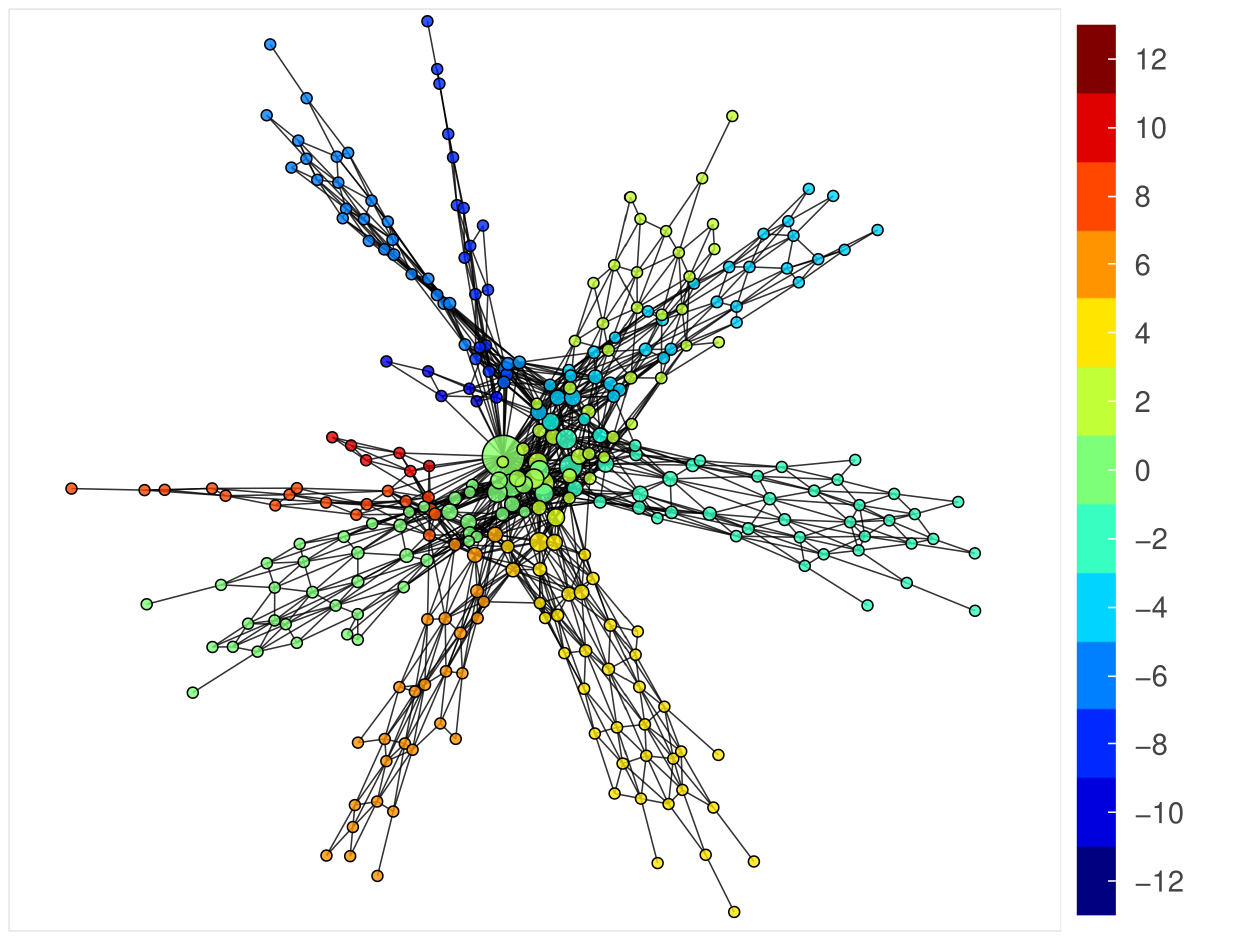}
		\caption{Khovanov for non-alternating}
		\label{fig:third}
	\end{subfigure}
	\caption{BallMapper plots for Khovanov homology for all knots up to 17 crossings (A), just non-alternating (C) and Jones data for not alternating (B) colored by the signature.}
	\label{fig:even_khov_BMSig}
\end{figure}

On the other hand,  Jones polynomial and signature of alternating knots  determine Khovanov homology (polynomial) \cite{garoufalidis2004conjecture}, while it is know that Khovanov homology is stronger invariant than the Jones polynomial e.g. it distinguishes knot $5_1$ and $10_{132}$ that have the same Jones polynomial \cite{KA}. 

From a theoretical point of view, the most intriguing  insight originates from \cite{dlotko2024mapper} that suggests that the flares of the Jones data Ball Mapper graph are related with the signature. More precisely, the flares, outside of the central clusters, contain only knots with equal signature. This insight relates with the Garoufalidis conjecture, building on the following theorem \cite{garoufalidis2003does, garoufalidis2005colored,garoufalidis2006non}.

\begin{thm}[Garoufalidis '03]
	For all simple knots (all the roots $\alpha\in\Delta(K)$  of the Alexander polynomial where $|\alpha|=1$  
	have multiplicity 1) up to 8 crossings and for all torus knots, the colored Jones polynomial determines the signature of the knot.
\end{thm}

\begin{figure}[h]
	\centering
	\begin{subfigure}{0.45\textwidth}
		\includegraphics[width=\textwidth]{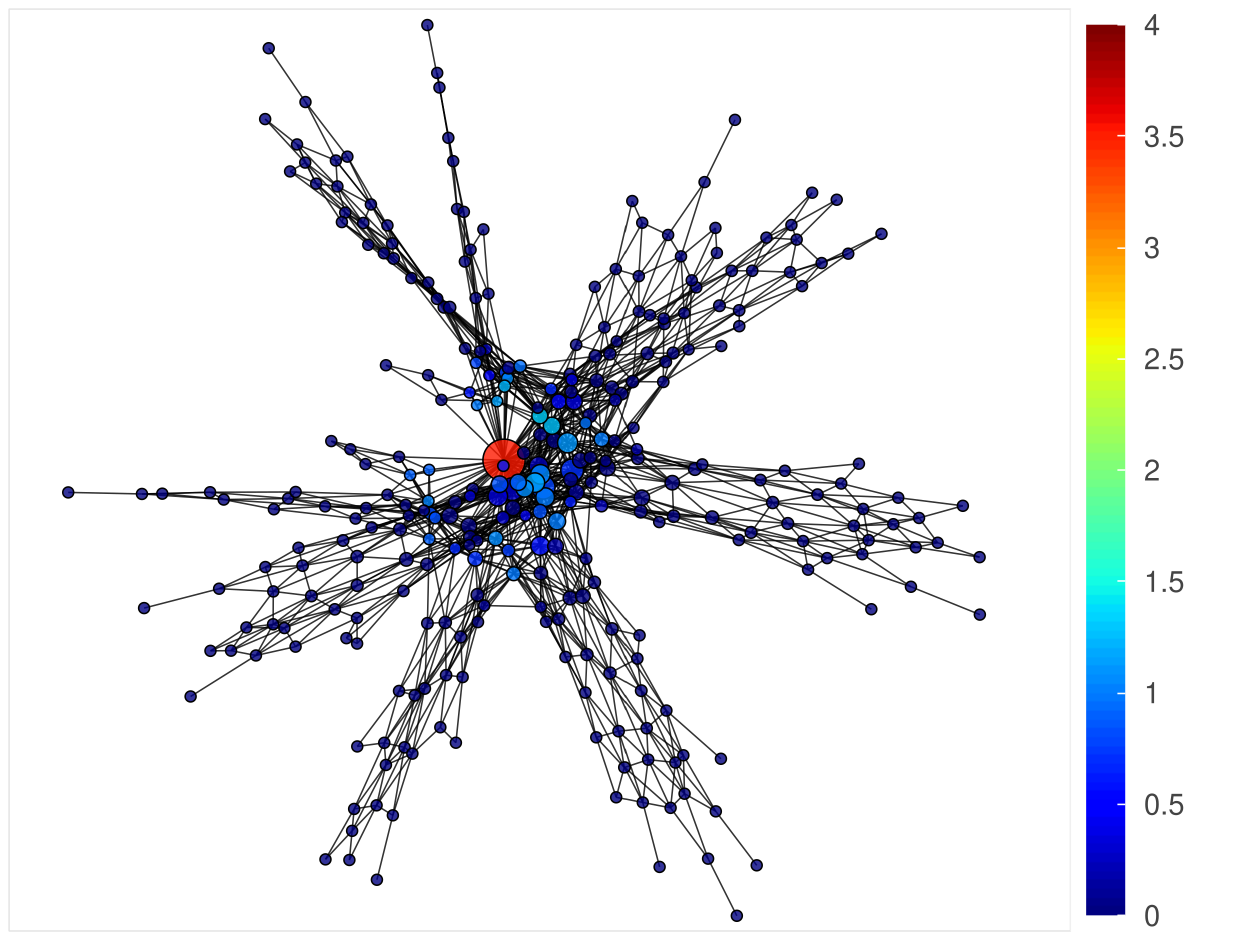}
		\caption{Khovanov.}
		\label{fig:first}
	\end{subfigure}
	\hfill
	\begin{subfigure}{0.45\textwidth}
		\includegraphics[width=\textwidth]{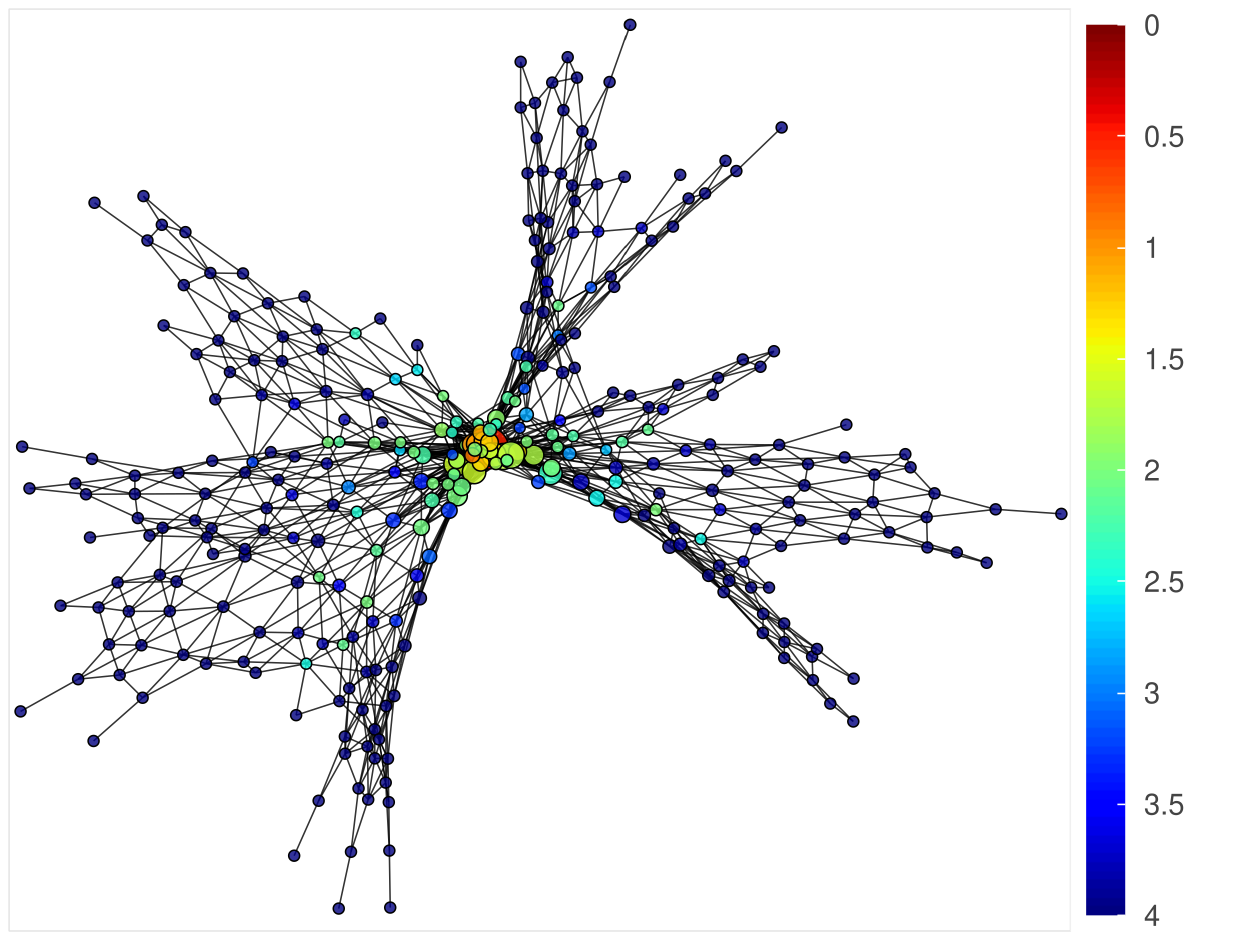}
		\caption{Jones.}
		\label{fig:JonesNA}
	\end{subfigure}
	
	\caption{BallMapper plots for even Khovanov homology (A) and Jones polynomial (B) for all non-alternating knots up to 17 crossings, colored by the standard deviation of signature in each node.}
	\label{fig:even_khov_BM_STD}
\end{figure}
Garoufalidis' conjecture was originally stated for the colored Jones polynomial and only  simple knots 
where simple is defined in terms of the Alexander polynomial.   
Therefore the observation that knots at the center of the Ball Mapper graph of a Jones polynomial happen to have mixed values of signature is not an evidence against Garoufalidis' conjecture, but an indication that perhaps,  in certain settings,  Jones polynomial relates with  signature, see Figures \ref{fig:BM_upto_15}, \ref{fig:even_khov_BMSig} showing Ball Mappers of Khovanov and Jones polynomial data colored by signature and the standard deviation of signature in each node Figure \ref{fig:even_khov_BM_STD}. 
For Khovanov-thin classes of knots where this would be true (up to some resolution determined by the Ball Mapper parameters) it would imply that the Jones polynomial alone determines the Khovanov homology (hence the polynomial) since it is proven that the Jones polynomial and signature determine Khovanov homology for all Khovanov thin knots  \cite{shumakovitch2004torsion}. However, in the light of the distribution of knots within the Ball Mapper graph, Ball Mapper insights into when does the Jones polynomial "determine" signature would have limited applicability.

\begin{figure}[h]
	\centering
	\begin{subfigure}{0.32\textwidth}
		\includegraphics[width=\textwidth]{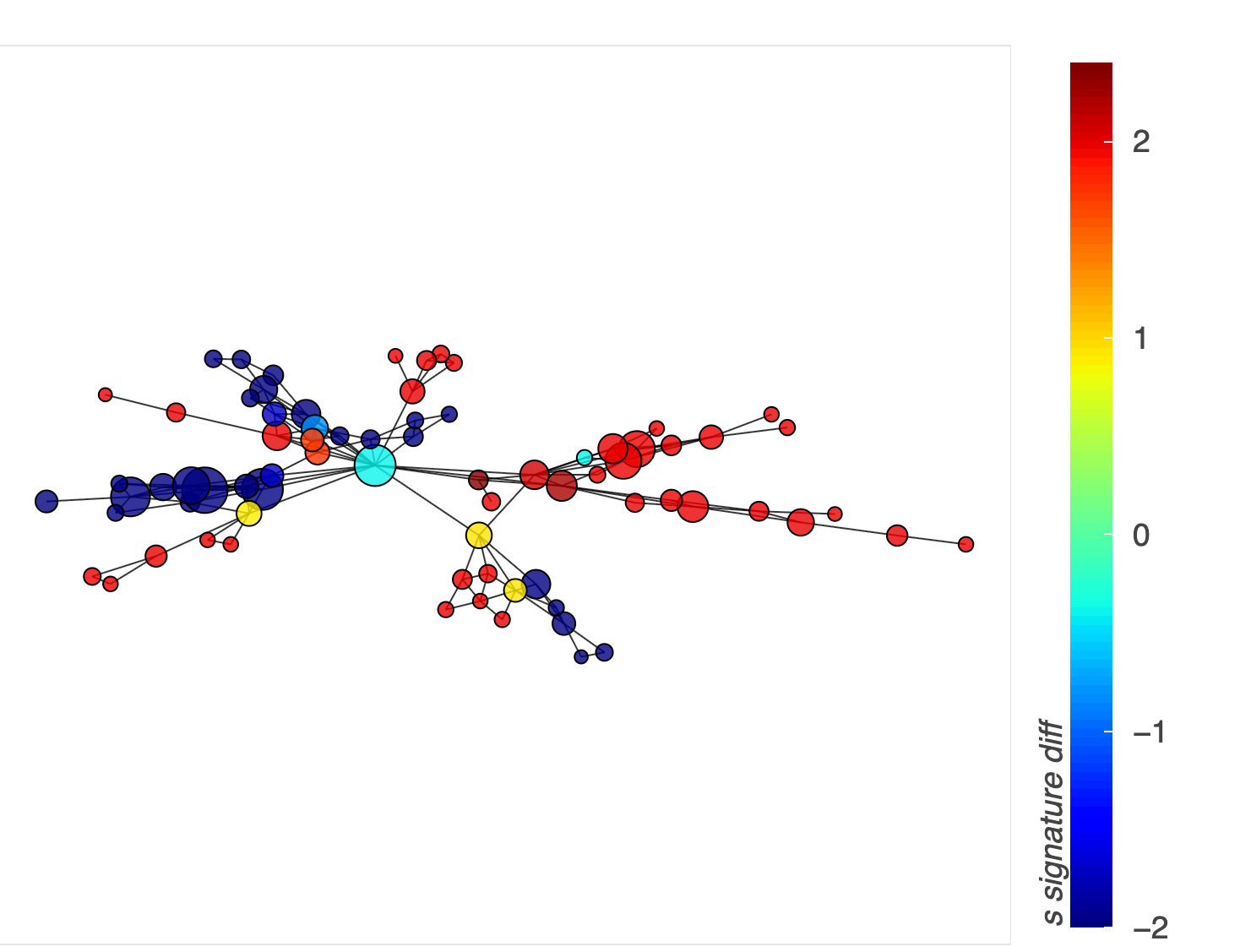}
		\caption{}
		\label{ssig}
	\end{subfigure}
	\hfill
	\begin{subfigure}{0.32\textwidth}
		\includegraphics[width=\textwidth]{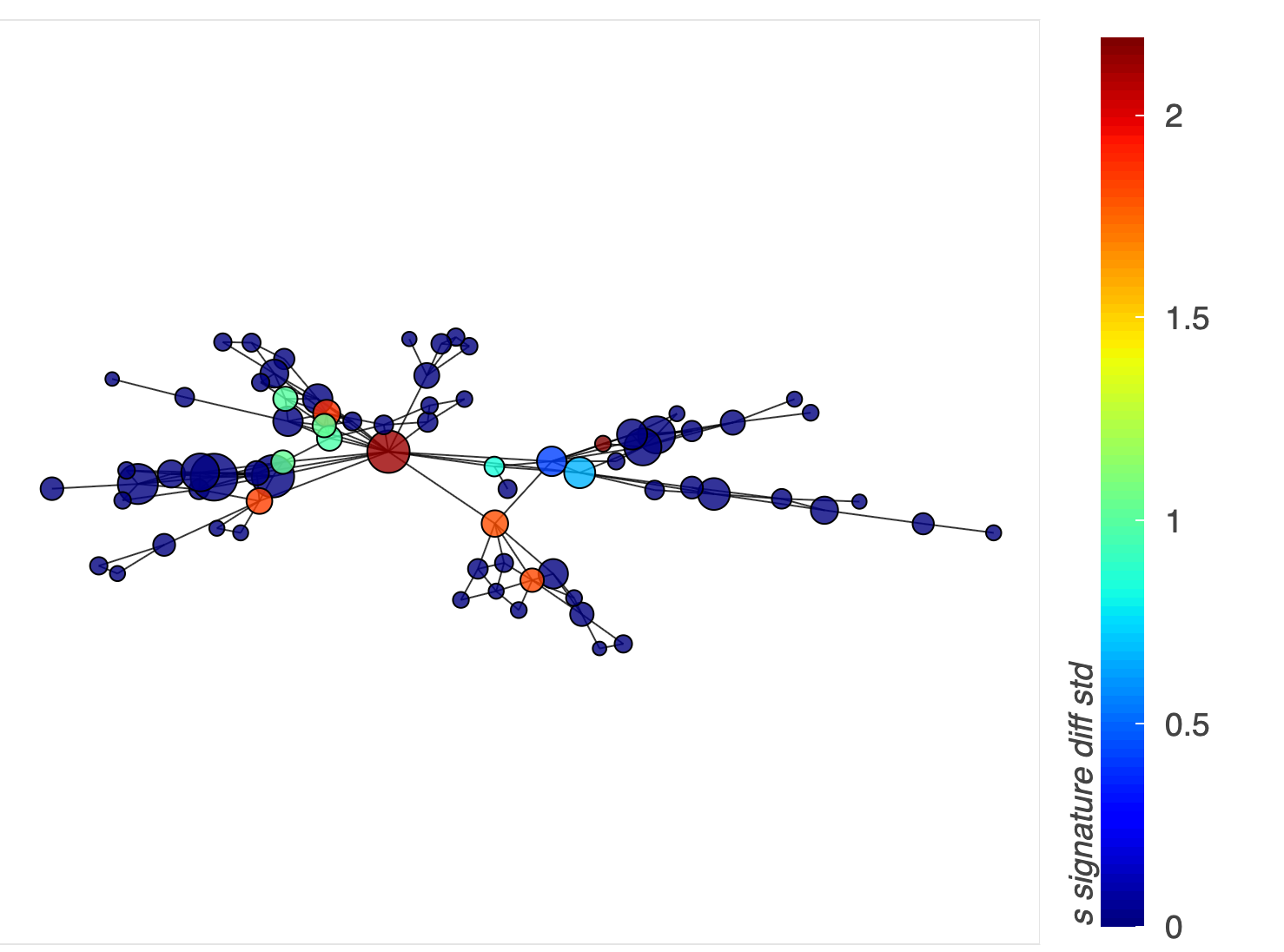}
		\caption{}
		\label{ssigstd}
	\end{subfigure}
	\hfill
	\begin{subfigure}{0.32\textwidth}
		\includegraphics[width=\textwidth]{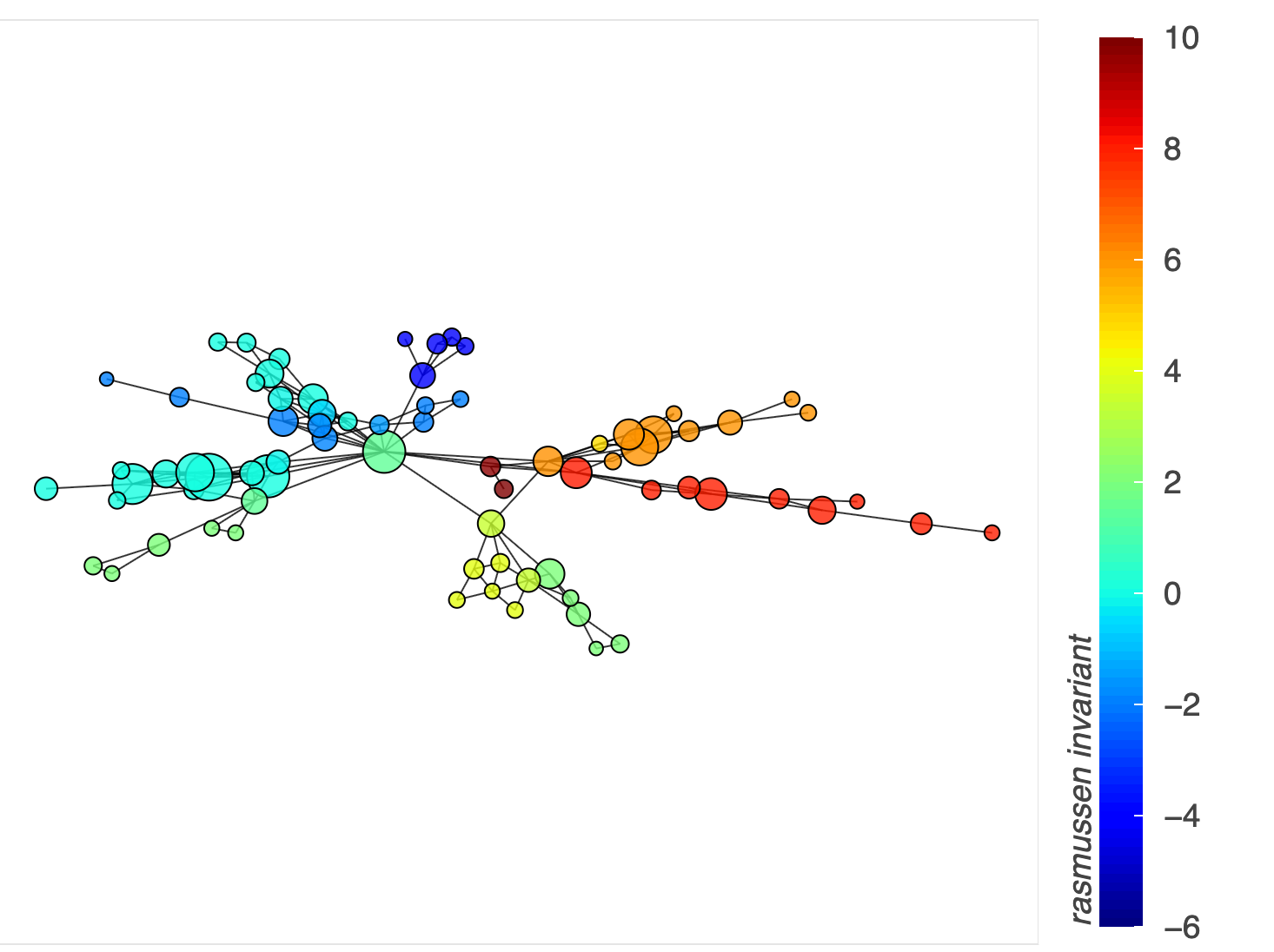}
		\caption{}\label{fig:rass}
	\end{subfigure}     
	\caption{Bal Mapper graphs for even Khovanov homology of 606 non-alternating knots up to 13 crossings where signature and s-invariant achieve different values colored by the difference between s and signature (A) and the average  standard deviation of the difference in each node (B), and Rasmussen s-invariant (C).}
	\label{fig:ssig}
\end{figure}

Recall that the Ball Mapper flares of Jones and Khovanov data correspond to sets of diagonals, signature and now we also have Rasmussen s-invariant Figure \ref{fig:rass}. Moreover, the difference is either 2 or 4 Figure \ref{ssig} and constant within flares away from the center Figure \ref{ssigstd}.

\begin{figure}[h]
	\begin{minipage}[b]{.49\linewidth}
		\centering
		\includegraphics[width=0.99\linewidth]{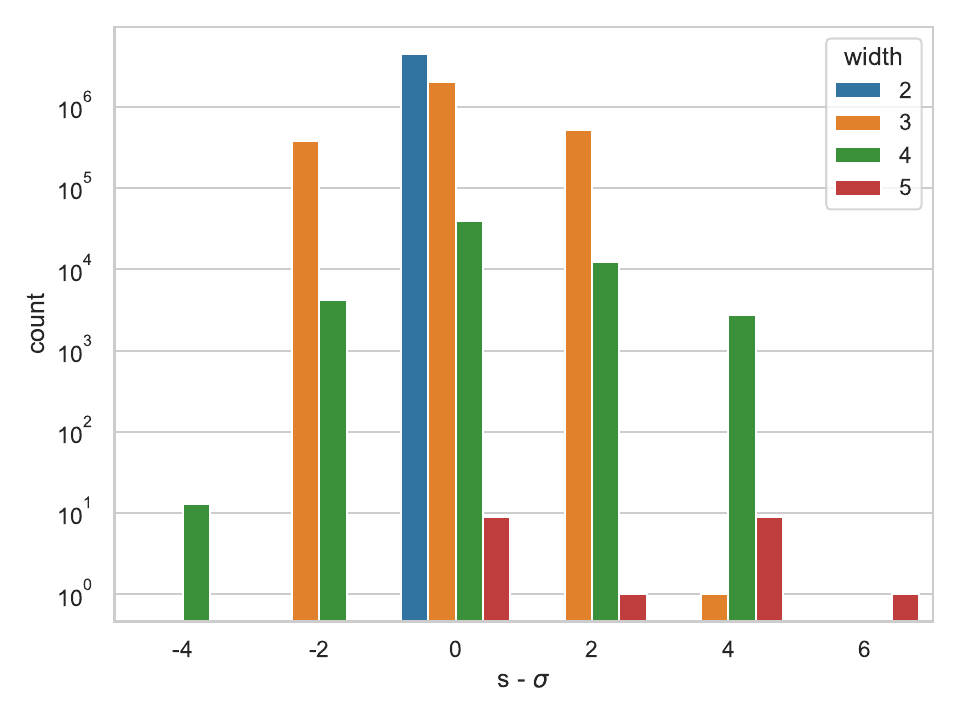}
		\captionof{figure}{The distribution of the difference between absolute s and absolute signature for all non-alternating knots up to 17 crossings, colored by the width of even Khovanov homology.}  \label{fig:s_vs_signature_width}
	\end{minipage}\hfill
	\begin{minipage}[b]{.49\linewidth}
		\centering
		\begin{tabular}{|c|c|c|c|}
			\hline
			& \textbf{$|s| = |\sigma|$} & \textbf{$|s| > |\sigma|$} & \textbf{$|s| < |\sigma|$} \\
			\hline
			8  & 3    & 0    & 0    \\
			9  & 10   & 0    & 1    \\
			10 & 45   & 6    & 2    \\
			11 & 216  & 11   & 11   \\
			12 & 999  & 80   & 47   \\
			13 & 5630 & 343  & 263  \\
			14 & 30099 & 2024 & 1549 \\
			15 & 179617 & 12674 & 9411 \\
			16 & 1070864 & 79643 & 60101 \\
			17 & 6572889 & 535553 & 385580 \\
			\hline
		\end{tabular}
		\captionof{table}{The number of non-alternating knots up various number of  crossings for which the absolute values of signature and s are the same or different.  }\label{kojebolji}
	\end{minipage}
\end{figure}

\begin{table}[!ht]
	\centering
	\begin{tabular}{|c|l|l|l|}
		\hline
		crossing no. & no. of  knots & no. of  knots $|s| \neq |\sigma|$ & percentage \\ \hline
		12 & 2978 & 127 & 4.26 \% \\ \hline
		13 & 12966 & 606 & 4.67 \% \\ \hline
		14 & 59938 & 3573 & 5.96 \% \\ \hline
		15 & 313231 & 22085 & 7.05 \% \\ \hline
		16 & 1701936 & 139744 & 8.21 \% \\ \hline
		17 & 9755329 & 921133 & 9.44 \% \\ \hline
	\end{tabular}
	\caption{The number of all knots up to various number of crossings and how many of them have $|s| \neq |\sigma|$. Note how the values in the third column corresponds to the sum of columns 2 and 3 of Table~\ref{kojebolji}.} \label{tab:sneSig}
\end{table}

Table \ref{tab:sneSig} provides data to support that the percentage of knots such that the absolute values of their signature and s-invariant coincide increases with the number of crossings and appears to be almost 10\% for 17 crossing knots.  

Figure \ref{fig:s_vs_signature_width} provides the distribution of the difference between s and signature for all knots up to 17 crossings but the lack of symmetry is at least partially due to the lack of mirrors. Knots with thin Khovanov homology in our data also have the same absolute values of signature and s, and higher difference relates with larger width of Khovanov homology.

\begin{conj}
	Absolute values of signature and s-invariant coincide for thin non-alternating knots. 
\end{conj}

As mentioned above, the conjecture is true for many classes of non-alternating knots \cite{manolescu2007khovanov} and we do not know how many of the non-alternating thin knots are not quasi-alternating.

Tables \ref{kojebolji} and \ref{tab:sneSig}  provide counts of knots with up to n crossings with the absolute value of the s-invariant  equal, bigger or smaller than  absolute value of signature. It shows that in the majority of the cases s and signature coincide and that they roughly outperform the other in providing the lower bound for the unknotting number with similar frequencies.

The following folklore conjecture has been around since the invention of Khovanov homology\footnote{RS thanks Alexander Shumakovitch and Adam Lowrance for bringing it to my attention.}
\begin{conj}\label{nase}
	If diagonals of integral Khovanov homology of $Kh(K, \mathbb{Z})$  are 
	$\{d_1, d_2, \ldots, d_w \}$, meaning that the width of its Khovanov homology is w,  then $d_1 <- \sigma(K) < d_w$. 
\end{conj}

Dasbach and Lowrance \cite{dasbach2011turaev} showed that a diagrammatic version of the conjecture holds for the Khovanov homology spanning tree complex. 

Data in Table \ref{tab:sigNotSpan} provides counterexamples to this Conjecture \ref{nase}. This Table includes the id of the knot, signature, and s, followed by the list of diagonals supporting  rational Khovanov homology $d(Kh, \mathbb{Q})$, diagonals featuring $Z_2$, then $Z_4$ torsion, followed by the width of Khovanov homology.  
Table \ref{tab:sigNotSpan} lists 4 knots up to 17 crossings that provide a counterexamples printed in boldface: 
$\mathbf{16n_{0229042}}$, 
$\mathbf{17nh_{0000038}}$, 
$\mathbf{17nh_{0000433}}$,  and $\mathbf{17nh_{0000273}}$,
as well as knots that provide counterexamples when considering only diagonals supporting the Khovanov homology with rational coefficients $(Kh, \mathbb{Q})$. 

\begin{table}[!ht]
	\centering
	\begin{tabular}{|l|l|l|l|l|l|l|}
		\hline
		Knot id & $\sigma$ & $s$ & $d(Kh, \mathbb{Q})$ & $Z_2$ diagonals & $Z_4$ diagonals & w  \\ \hline
		$15n_{041127}$ & 0 & 2  &\{1, 3, 5\}  & \{1, 3, -1 \} &   & 4 
		%& 3
		\\ \hline
		$16n_{0196771}$ & 6 & 10  & \{9, 11, 7  & \{9, 5, 7 \} & \{7\} & 4 
		%& 2
		\\ \hline
		$\mathbf{16n_{0229042}}$ & 0 & -2  & \{-5, -3, -1\}  & \{-5, -3, -1 \} &   & 3 %& 3 
		\\ \hline
		$16n_{0245346}$ & -2 & 0  & \{1, 3, -1 \} & \{1, -3, -1 \} &   & \{4\} %& 3 
		\\ \hline
		$\mathbf{17nh_{0000038}}$ & 0 & -2  & \{-5, -3, -1 \} & \{-5, -3, -1\}  &   & 3
		%& 2 
		\\ \hline
		$17nh_{0000075}$ & 8 & 12  & \{9, 11, 13 \} & \{9, 11, 7 \} & 9  & \{4 \}
		%& 2 
		\\ \hline
		$\mathbf{17nh_{0000124}}$ & 2 & 0  & \{1, -3, -1]  & \{1, -3, -1 ] &   & 3 
		%& 2 
		\\ \hline
		$\mathbf{17nh_{0000158}}$ & -8 & -12  & \{-13, -11, -9]  &\{ -13, -11, -9 ] & \{-11\} & 3 
		%& 2
		\\ \hline
		$17nh_{0000173}$ & 0 & 2  & \{1, 3, 5 \} & \{1, 3, -1 \} &   & 4
		%& 3 
		\\ \hline
		$17nh_{0000192}$ & 0 & 2  & \{1, 3, 5 \} & \{1, 3, -1\}  &   & 4 %& 3 
		\\ \hline
		$\mathbf{17nh_{0000273}}$ & 0 & -2  & \{-5, -3, -1\}  & \{-5, -3, -1\}  &   & 3 %& 3 
		\\ \hline
		$17nh_{0000281}$ & -2 & 0  & \{1, 3, -1 \} & \{1, -3, -1 \} &   & 4 %& 3
		\\ \hline
		$17nh_{0000283}$ & 6 & 10  &\{ 9, 11, 7\}  & \{9, 5, 7 \} & \{7\} & 4 %& 3
		\\ \hline
		$17nh_{0000297}$ & 0 & 2  & \{1, 3, 5 \} & \{1, 3, -1 \} &   & 4 %& 3 
		\\ \hline
		$17nh_{0000410}$ & 0 & 2  & \{1, 3, 5\}  & 1, 3, -1 \} &   & 4 %& 3 
		\\ \hline
		$\mathbf{17nh_{0000433}}$ & 0 & -2  & \{-5, -3, -1 \} & \{-5, -3, -1\}  &   & 3 %& 4 
		\\ \hline
	\end{tabular}
	
	\caption{Knots with at most 17 crossings whose signature does not fall in the span of the diagonals for the free part of Khovanov homology yet they are on adjacent diagonals. Bold font denotes knots with the same property for integral Khovanov homology.}\label{tab:sigNotSpan}
\end{table}

\begin{table}[ht]
	\centering
	\begin{tabular}{|l|l|l|l|l|l|l|l|}
		\hline
		Knot id & $\sigma$ & $s$ & $d(Kh, \mathbb{Q})$ & $Z_2$ diagonals & $Z_4$ diagonals & $d(Kh, \mathbb{Z})$ &w  \\ \hline
		$17nh_{0000460}$ & 8 & 14 & \{7, 9, 11, 13, 15\} & \{9, 11, 13, 7\} & 9 & \{7, 9, 11, 13, 15\} & 5   \\ \hline
	\end{tabular}
	\caption{The only knot up to 17 crossings with  difference between s and signature is equal to 6.}\label{theonly6}
\end{table}

Notice that there is only one knot up to 17 crossings where the difference between absolute values of s and signature is 6, is $17nh_{0000460}$ with the invariants listed in Table \ref{theonly6}. It's Khovanov homology has width 5 realized by diagonals $\{ 7, 9, 11, 13, 15\}$ that supports value of $s=14$ and $\sigma=8$ and in zeroth homological grading it is supported in gradings 13 and 15. 

%\newpage
%%%%%%%%%%%%%%%%%%%%%%%%%%%%%
\section{Knot data from random polygons}\label{random} 

\begin{figure}[h!]
	\centering
	\begin{subfigure}[b]{0.45\textwidth}
		\includegraphics[width=\textwidth]{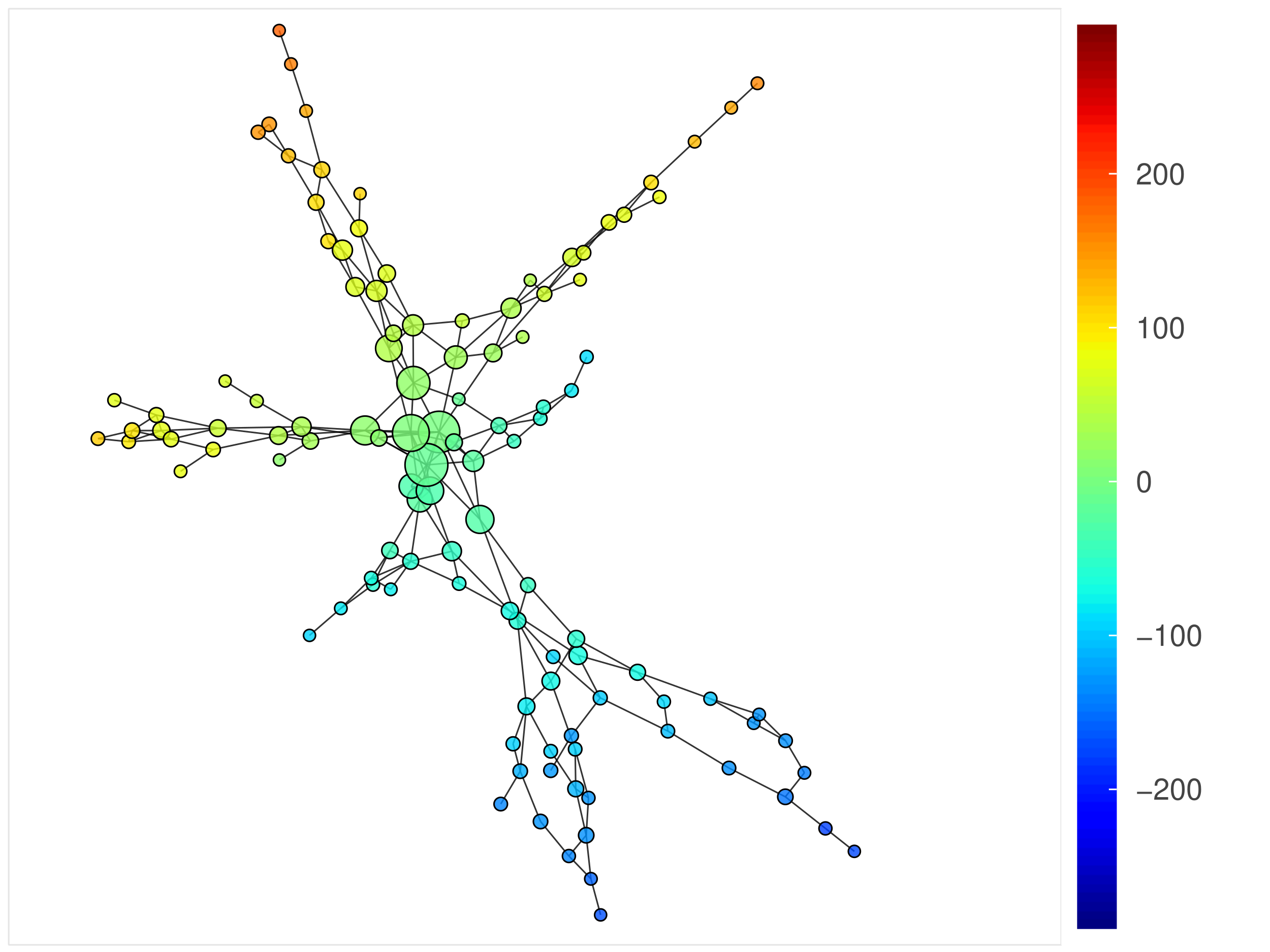}
		\caption{}\label{randomdet}
	\end{subfigure}
	\hfill
	\begin{subfigure}[b]{0.45\textwidth}
		\includegraphics[width=\textwidth]{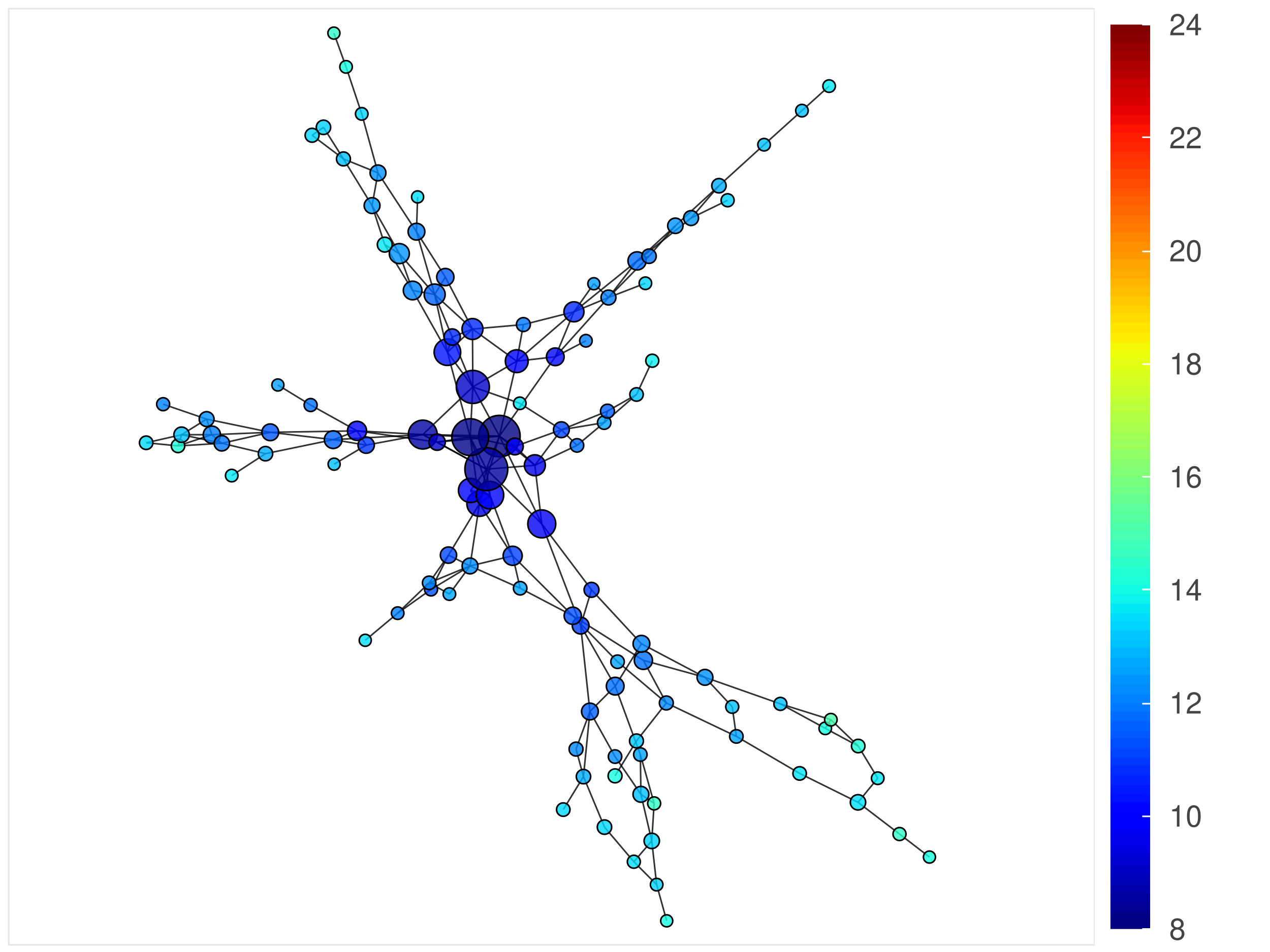}
		\caption{}\label{randomspan}
	\end{subfigure}
	\caption{Ball Mapper graphs of the Jones polynomial data of random knots colored by (A) determinant and (B) span of the Jones polynomial. Note that knots with extremely high span of the Jones polynomial and determinant values are disconnected and omitted.}
	
	\label{fig:jones_random}
\end{figure}

In order to confirm and potentially generalize observations obtained using topological exploratory and visualization tools we construct Ball Mapper on the Jones polynomial data of about 10.000 knots obtained from random polygons of length 500 as described in Subsection \ref{randomdata}. This s data set contained  2572 of distinct polynomials., which is about 26\%, similar to the percentage of unique Jones polynomials reported  in Table \ref{tab:DistUni}.  Ball Mapper of Jones polynomial for random knots shares the same star-like structure after many disconnected components are removed Figure \ref{randomdet}.
We color these Ball Mapper graph by the determinant Figure\ref{randomdet} and span of the Jones polynomial Figure\ref{randomspan} to conform our observations that the absolute value of the determinant grows when we move radially away from the center. However, the mean value of the determinant was about  922320, with huge standard deviation since the maximal value of the determinant on this data set was  1264326765! Similarly the average span was 17.28 ranging all the way up to 53! Such a high value confirms that the data contained at least some complicated knots. These results warrant a thorough investigation of a larger random knots data set, especially in the light of recent results on knot recognition based on geometric writhe \cite{sleiman2024geometric}, and relations between writhe, signature and ropelength \cite{klotz2025ropelength, thompson2025space}.

\section*{Acknowledgments}
The authors want to thank Dirk Shuetz and Alexander Shumakovitch for providing the data for the Rasmussen s-invariant and Khovanov homology used in this paper and Adam Lowrance for insightful comments on the analysis of knot--theoretic data. 
PD and DG acknowledge support by Dioscuri program initiated by the Max Planck Society, jointly managed with the National Science Centre (Poland), and mutually funded by the Polish Ministry of Science and Higher Education and the German Federal Ministry of Education and Research. DG is also with the University of Warsaw within the Doctoral School of Exact and Natural Sciences. RS is partially supported by NSF grant DMS-1854705. The authors declare no competing interests.

\end{document}